\documentclass[12pt]{amsart}
\usepackage{amsmath, amssymb, latexsym, amsthm, mathrsfs, color, mathtools} 
\usepackage[all]{xy}
\usepackage[top=2in, bottom=1.5in, left=1in, right=1in]{geometry}

\newenvironment{changemargin}[2]{\begin{list}{}{
\setlength{\topsep}{0pt}\setlength{\leftmargin}{0pt}
\setlength{\rightmargin}{0pt}
\setlength{\listparindent}{\parindent}
\setlength{\itemindent}{\parindent}
\setlength{\parsep}{0pt plus 1pt}
\addtolength{\leftmargin}{#1}\addtolength{\rightmargin}{#2}
}\item }{\end{list}}

\newcommand{\beq}{\begin{eqnarray*}}\newcommand{\eeq}{\end{eqnarray*}}
\newcommand{\mbq}{\mb{?}}
\newcommand{\mb}[1]{{\mbox{\textbf{#1}}}}
\newcommand{\nop}{$\times$}
\newcommand{\fbn}{\!\!\fbox{\!\nop\!}\!\!}
\newcommand{\yup}{\checkmark}
\newcommand{\forces}{\Vdash}
\newcommand{\name}[1]{\dot{#1}}

\newcommand{\FU}{Fr\'echet--Urysohn}
\newcommand{\gs}{$\gamma$~space}
\newcommand{\Ga}{\Gamma}\newcommand{\Om}{\Omega}
\newcommand{\smallbinom}[2]{\begin{psmallmatrix} #1\\ #2 \end{psmallmatrix}}
\newcommand{\bgamma}{\smallbinom{\Om}{\Ga}}
\newcommand{\two}{\{0,1\}}

\newcommand{\Sel}{\mathsf{S}}
\newcommand{\sset}[2]{{\{\,#1 : #2\,\}}}
\newcommand{\smb}[1]{{\!\!\mb{#1}\!\!}}
\newcommand{\medset}[2]{{\biggl\{\,#1 : #2\,\biggr\}}}
\newcommand{\smallmedset}[2]{{\bigl\{\,#1 : #2\,\bigr\}}}
\newcommand{\set}[2]{{\left\{\,#1 : #2\,\right\}}}
\newcommand{\seq}[2]{{\la\, #1 : #2\,\ra}}

\newcommand{\concat}[1]{\hat{\phantom{a}}\langle #1\rangle}
\newcommand{\poset}{\mathbb{P}}

\newcommand{\alephes}{{\aleph_0}}

\newcommand{\Cp}{\op{C}_\mathrm{p}}

\newcommand{\Pa}[8]{\bibitem{#1} {#2}, \emph{#3}, {#4} \textbf{#5} ({#6}), {#7}--{#8}.}

\newcommand{\Bc}[9]{\bibitem{#1} {#2}, \emph{#3}, in: \textbf{#4} (#5), #6 #7, #8--#9.}

\newcommand{\od}{\mathfrak{od}}
\newcommand{\Setting}[7]{\xymatrix@R=4pt@C=7pt{#1\ar@{-}[r]&#2\ar@{-}[r]&#3\\&#4\ar@{-}[u]\\
#5\ar@{-}[uu]\ar@{-}[r] & #6\ar@{-}[u]\ar@{-}[r] & #7\ar@{-}[uu]}}

\newcommand{\bq}{\begin{quote}}
\newcommand{\eq}{\end{quote}}
\newcommand{\cl}[1]{\overline{#1}}
\newcommand{\CH}{the Continuum Hypothesis}

\newcommand{\inv}{^{-1}}
\newcommand{\Cantor}{{\two^\N}}

\newcommand{\sr}[3]{\underset{\mbox{#3}}{\mbox{#1}}}
\newcommand{\gp}{\binom{\Om}{\Ga}}

\newcommand{\N}{\mathbb{N}}
\newcommand{\NN}{{\N^{\N}}}

\newcommand{\PN}{{P(\N)}}
\newcommand{\roth}{{[\N]^{\mbox{\tiny $\infty$}}}} 
\newcommand{\Fin}{[\N]^{\mbox{\tiny $<\!\infty$}}} 

\newcommand{\sseq}[1]{\{#1 : n\in\N\}}

\newcommand{\op}{\operatorname}

\newcommand{\scrA}{\mathscr{A}}
\newcommand{\scrB}{\mathscr{B}}

\newcommand{\B}{\mathrm{B}}
\newcommand{\cB}{\mathcal{B}}

\newcommand{\BG}{\B_\Ga}

\newcommand{\BT}{\B_\Tau}

\newcommand{\BO}{\B_\Om}

\newcommand{\Tau}{\mathrm{T}}
\newcommand{\cA}{\mathcal{A}}

\newcommand{\cD}{\mathcal{D}}
\newcommand{\cF}{\mathcal{F}}

\newcommand{\cM}{\mathcal{M}}

\newcommand{\cO}{\mathcal{O}}
\newcommand{\Op}{\mathrm{O}}
\newcommand{\rmA}{\mathrm{A}}
\newcommand{\rmB}{\mathrm{B}}
\newcommand{\rmD}{\mathrm{D}}

\newcommand{\R}{\mathbb{R}}
\newcommand{\cU}{\mathcal{U}}
\newcommand{\Un}{\bigcup}
\newcommand{\cV}{\mathcal{V}}
\newcommand{\cW}{\mathcal{W}}

\newcommand{\Impl}{\Rightarrow}
\long\def\forget#1\forgotten{}
\newcommand{\ft}{\mathfrak{t}}
\newcommand{\fb}{\mathfrak{b}}

\newcommand{\fd}{\mathfrak{d}}

\newcommand{\oo}{\infty}

\newcommand{\fp}{\mathfrak{p}}

\newcommand{\fs}{\mathfrak{s}}

\newcommand{\x}{\times}

\newcommand{\comp}{^{\text{\tt c}}}

\newcommand{\sub}{\subseteq}
\newcommand{\spst}{\supseteq}
\newcommand{\sm}{\setminus}
\newcommand{\as}{\subseteq^*}

\newcommand{\la}{\langle}
\newcommand{\ra}{\rangle}

\newcommand{\cov}{\op{cov}}

\newcommand{\non}{\op{non}}

\newtheorem{thm}{Theorem}[section]
\newcommand{\bthm}{\begin{thm}} \newcommand{\ethm}{\end{thm}}
\newtheorem{prop}[thm]{Proposition}
\newcommand{\bprp}{\begin{prop}} \newcommand{\eprp}{\end{prop}}
\newtheorem{fact}[thm]{Fact}
\newcommand{\bfct}{\begin{fact}} \newcommand{\efct}{\end{fact}}
\newtheorem{prob}[thm]{Problem}
\newcommand{\bprb}{\begin{prob}} \newcommand{\eprb}{\end{prob}}
\newtheorem{lem}[thm]{Lemma}
\newcommand{\blem}{\begin{lem}} \newcommand{\elem}{\end{lem}}
\newtheorem{claim}[thm]{Claim}
\newcommand{\bclm}{\begin{claim}} \newcommand{\eclm}{\end{claim}}
\newtheorem{cor}[thm]{Corollary}
\newcommand{\bcor}{\begin{cor}} \newcommand{\ecor}{\end{cor}}
\newtheorem{conj}[thm]{Conjecture}
\newcommand{\bcnj}{\begin{conj}} \newcommand{\ecnj}{\end{conj}}
\theoremstyle{definition}
\newtheorem{defn}[thm]{Definition}
\newcommand{\bdfn}{\begin{defn}} \newcommand{\edfn}{\end{defn}}
\theoremstyle{remark}
\newtheorem{rem}[thm]{Remark}
\newcommand{\brem}{\begin{rem}} \newcommand{\erem}{\end{rem}}
\newtheorem{cnv}[thm]{Convention}
\newcommand{\bcnv}{\begin{cnv}} \newcommand{\ecnv}{\end{cnv}}
\newtheorem{exam}[thm]{Example}
\newcommand{\bexm}{\begin{exam}} \newcommand{\eexm}{\end{exam}}
\newcommand{\bpf}{\begin{proof}} \newcommand{\epf}{\end{proof}}
\newcommand{\be}{\begin{enumerate}}
\newcommand{\ee}{\end{enumerate}}
\newcommand{\bi}{\begin{itemize}}
\newcommand{\itm}{\item}
\newcommand{\ei}{\end{itemize}}

\forget
\forgotten

\newcommand{\sone}{\mathsf{S}_1}
\newcommand{\sfin}{\mathsf{S}_\mathrm{fin}}
\newcommand{\ufin}{\mathsf{U}_\mathrm{fin}}

\newcommand{\gone}{\mathsf{G}_1}

\newcommand{\ed}{

\subsection*{Acknowledgments}
The research of the third named author was supported by FWF Grant  I 1209-N25 and
the program APART of the Austrian Academy of Sciences.
Parts of the work reported here were carried out during three visits of the second named author
at the Kurt G\"odel Research Center for Mathematical Logic.
These visits were partially supported by the above-mentioned grants and by
the Research Networking Programme \emph{New Frontiers of Infinity} (INFTY),
funded by the European Science Foundation.
The second named author thanks the third named author for his kind hospitality,
and the Kurt G\"odel Research Center Director, researchers and staff for the excellent
academic and friendly atmosphere.
We thank Francis Jordan for useful discussions, and the referee for his/her work on
reviewing this paper.

\end{document}}

\title[Products of $\gamma$~spaces]{Selective covering properties of product spaces, II:\\
$\gamma$~spaces}
\author[A. Miller]{Arnold W. Miller}
\address[Miller]{Department of Mathematics, University of Wisconsin-Madison, Van Vleck Hall 480 Lincoln Drive, Madison, Wisconsin 53706-1388, USA}
\email{miller@math.wisc.edu}
\urladdr{http://www.math.wisc.edu/~miller}
\author[B. Tsaban]{Boaz Tsaban}
\address[Tsaban]{Department of Mathematics, Bar-Ilan University, Ramat Gan 5290002, Israel, and
Faculty of Mathematics and Computer Science, Weizmann Institute of Science, Rehovot 7610001, Israel}
\email{tsaban@math.biu.ac.il}
\urladdr{http://math.biu.ac.il/~tsaban}
\author[L. Zdomskyy]{Lyubomyr Zdomskyy}
\address[Zdomskyy]{Kurt G\"odel Research Center for Mathematical Logic, University of Vienna,
W\"ahringer Str.\ 25, 1090 Vienna, Austria}
\email{lzdomsky@logic.univie.ac.at}
\urladdr{http://www.logic.univie.ac.at/~lzdomsky}



\begin{document}


\begin{abstract}
We study productive properties of \gs{}s, and their relation to other, classic and modern, selective covering properties.
Among other things, we prove the following results:
\be
\item Solving a problem of F. Jordan, we show that for every unbounded tower set $X\sub\R$ of cardinality $\aleph_1$,
the space $\Cp(X)$ is productively \FU{}. In particular, the set $X$ is productively $\gamma$.
\item Solving problems of Scheepers and Weiss, and proving a conjecture of Babinkostova--Scheepers, we prove that,
assuming \CH{}, there are \gs{}s whose product is not even Menger.
\item Solving a problem of Scheepers--Tall, we show that
the properties $\gamma$ and Gerlits--Nagy (*) are preserved by Cohen forcing. Moreover, every Hurewicz space that remains
Hurewicz in a Cohen extension must be Rothberger (and thus (*)).
\ee
We apply our results to solve a large number of additional problems,
and use Arhangel'ski\u{\i} duality to obtain results concerning local properties of function spaces and countable topological
groups.
\end{abstract}

\maketitle

\section{Introduction}

For a Tychonoff space $X$, let $\Cp(X)$ be the space of continuous real-valued functions on $X$,
endowed with the topology of pointwise convergence, that is, the topology inherited from the
Tychonoff product $\R^X$.
In their seminal paper~\cite{GN}, Gerlits and Nagy characterized the property that the space $\Cp(X)$ is \emph{\FU{}}---that every
point in the closure of a set is the limit of a sequence of elements from that set---in terms of
a covering property of the domain space $X$.
We study the behavior of this covering property under taking products with spaces possessing related
covering properties.

By \emph{space} we mean an infinite topological space.
Whenever the space $\Cp(X)$ is considered, we tacitly restrict our scope to Tychonoff
spaces.
The concrete examples constructed in this paper are all subsets of the real line.

The covering property introduced by Gerlits and Nagy is best viewed in terms of its relation
to other, selective covering properties.
The framework of \emph{selection principles} was introduced by Scheepers in~\cite{coc1}
to study, in a uniform manner, a variety of properties introduced
in several mathematical contexts, since the early 1920's.
Detailed introductions are available in~\cite{KocSurv, LecceSurvey, ict, SaSchRP}.
We provide here a brief one, adapted from~\cite{LinSAdd}.

Let $X$ be a space. We say that $\cU$ is a \emph{cover}
of $X$ if $X=\Un\cU$, but $X$ is not covered by any single member of $\cU$.
Let $\Op(X)$ be the family of all open covers of $X$.
When $X$ is considered as a subspace of a larger space $Y$,
the family $\Op(X)$ consists of the covers of $X$ by open subsets of $Y$.
Define the following subfamilies of $\Op(X)$:
$\cU\in\Om(X)$ if each finite subset of $X$ is contained in some member of $\cU$.
$\cU\in\Ga(X)$ if $\cU$ is infinite, and each element of $X$ is contained in all but
finitely many members of $\cU$.

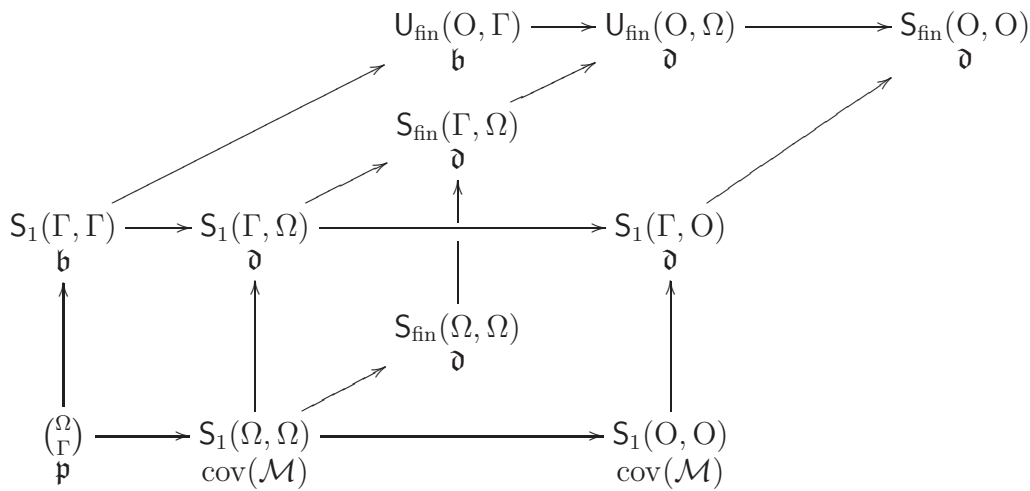
\begin{figure}[!htp]
\begin{changemargin}{-4cm}{-3cm}
\begin{center}
{
$\xymatrix@R=8pt{
&
&
& \sr{$\ufin(\Op,\Ga)$}{Hurewicz}{$\fb$}\ar[r]
& \sr{$\ufin(\Op,\Om)$}{}{$\fd$}\ar[rr]
& & \sr{$\sfin(\Op,\Op)$}{Menger}{$\fd$}
\\
&
&
& \sr{$\sfin(\Ga,\Om)$}{}{$\fd$}\ar[ur]
\\
& \sr{$\sone(\Ga,\Ga)$}{}{$\fb$}\ar[r]\ar[uurr]
& \sr{$\sone(\Ga,\Om)$}{}{$\fd$}\ar[rr]\ar[ur]
& & \sr{$\sone(\Ga,\Op)$}{}{$\fd$}\ar[uurr]
\\
&
&
& \sr{$\sfin(\Om,\Om)$}{}{$\fd$}\ar'[u][uu]
\\
& \sr{$\gp$}{Gerlits--Nagy}{$\fp$}\ar[r]\ar[uu]
& \sr{$\sone(\Om,\Om)$}{}{$\cov(\cM)$}\ar[uu]\ar[rr]\ar[ur]
& & \sr{$\sone(\Op,\Op)$}{Rothberger}{$\cov(\cM)$}\ar[uu]
}$
}
\caption{The Scheepers Diagram}\label{SchDiag}
\end{center}
\end{changemargin}
\end{figure}

Some of the following statements may hold for families $\scrA$ and $\scrB$ of covers of $X$.
\begin{description}
\item[$\smallbinom{\scrA}{\scrB}$] Each member of $\scrA$ contains a member of $\scrB$.
\item[$\sone(\scrA,\scrB)$] For each sequence $\seq{\cU_n\in\scrA}{n\in\N}$, there is a selection
$\seq{U_n\in\cU_n}{n\in\N}$ such that $\sseq{U_n}\in\scrB$.
\item[$\sfin(\scrA,\scrB)$] For each sequence $\seq{\cU_n\in\scrA}{n\in\N}$, there is a selection
of finite sets $\seq{\cF_n\sub\cU_n}{n\in\N}$ such that $\Un_n\cF_n\in\scrB$.
\item[$\ufin(\scrA,\scrB)$] For each sequence $\seq{\cU_n\in\scrA}{n\in\N}$, where no $\cU_n$ contains a finite
subcover, there is a selection
of finite sets $\seq{\cF_n\sub\cU_n}{n\in\N}$ such that $\sseq{\Un\cF_n}\in\scrB$.
\end{description}
We say, e.g., that $X$ satisfies $\sone(\Op,\Op)$ if the statement $\sone(\Op(X),\Op(X))$ holds.
This way, the notation $\sone(\Op,\Op)$ stands for a property (or a class) of spaces.
An analogous convention is followed for all other selection principles and families of covers.
Each nontrivial property among these properties, where $\scrA$ and $\scrB$ range over $\Op,\Om$ and $\Ga$,
is equivalent to one in Figure~~\ref{SchDiag}~\cite{coc1, coc2}. Some of the equivalences request
that the space be Lindel\"of. All spaces constructed in this paper to satisfy properties in the Scheepers Diagram are Lindel\"of.
Moreover, they are all subspaces of $\R$.

In the Scheepers Diagram, an arrow denotes implication.
We indicate below each class $P$ its \emph{critical cardinality} $\non(P)$,
the minimal cardinality of a space not in the class.
These cardinals are all combinatorial cardinal characteristics of the continuum,
details about which are available in~\cite{BlassHBK}. Following the convention in the
field of selection principles, influenced by the monograph~\cite{BarJu},
we deviate from the notation in~\cite{BlassHBK} by denoting the family of meager (Baire first category) 
sets in $\R$ by $\cM$.

The properties $\ufin(\Op,\Ga)$, $\sfin(\Op,\Op)$ and $\sone(\Op,\Op)$ were first studied by
Hurewicz, Menger and Rothberger, respectively.
\emph{\gs{}s} were introduced by Gerlits and Nagy~\cite{GN} as the spaces satisfying $\bgamma$.
Gerlits and Nagy proved that, for a space $X$, the space $\Cp(X)$ is \FU{}
if and only if $X$ is a \gs{}.

We also consider the classes of covers
$\B$, $\BO$ and $\BG$, defined as $\Op$, $\Om$ and $\Ga$ were
defined, replacing \emph{open cover} by \emph{countable Borel cover}.
The properties thus obtained have rich history of their own~\cite{CBC}, and
the Borel variants of the studied properties are strictly stronger than the open ones~\cite{CBC}.
Many additional---classic and new---properties are studied in relation to the Scheepers Diagram.

\bdfn
Let $P$ be a property (or class) of spaces.
A space $X$ is \emph{productively $P$} if $X\x Y$ has the property $P$ for each space $Y$ satisfying $P$.
\edfn


In Section~\ref{sec:prod} we construct productively \gs{}s in $\R$ from a weak hypothesis.
In Section~\ref{sec:nonprod} we construct, using \CH{}, two \gs{}s in $\R$ whose product is not Menger.
In Section~\ref{sec:apps} we use our results to solve a large number of problems, from the literature and
from the folklore of selection principles.
In Section~\ref{sec:forcing} we determine the effect of Cohen forcing on \gs{}s, Hurewicz spaces, and Gerlits--Nagy (*) spaces.
In Section~\ref{sec:cp} we use our results together with $\Cp$ theory, to obtain new results concerning local and density
properties of function spaces. In the last section, we prove that every product of an unbounded tower set and a
Sierpi\'nski set satisfies $\sone(\Ga,\Ga)$.

\section{Productively \gs{}s in $\R$}\label{sec:prod}

Recall the Gerlits--Nagy Theorem that a space $X$ is a \gs{} if and only if the space $\Cp(X)$ is \FU{}.
In his papers~\cite{Jor07, Jor08}, F. Jordan studied the property that $\Cp(X)$ is \emph{productively}
\FU{}. (In this case, it is said in~\cite{Jor07, Jor08} that the space
$X$ is a \emph{productive $\gamma$-space}.
Since this terminology is admitted in~\cite{Jor08} to be confusing, we avoid it here.)
We begin with a short proof of a result of Jordan.
In the proof, and later on,
we use the following observations.

\blem\label{lem:un2prod}
Let $P$ be a class of spaces that is hereditary for closed subsets and is preserved by finite powers. Then
for all spaces $X$ and $Y$ such that the disjoint union space $X\sqcup Y$ satisfies $P$,
the product space $X\x Y$ satisfies $P$, too. In particular,
if $P$ is preserved by finite unions, then it is preserved by finite products.
\elem
\bpf
We prove the first assertion.
As $P$ is preserved by finite powers, the space $(X\sqcup Y)^2$ satisfies $P$.
As $X\x Y$ is a closed subset of $(X\sqcup Y)^2$, it satisfies $P$, too.
\epf

If the disjoint union space $X\sqcup Y$ is a \gs{}, then so is the product space $X\x Y$~\cite[Proposition 2.3]{Miller07}.

\bcor\label{cor:unprod}
Let $X$ and $Y$ be spaces. The disjoint union space $X\sqcup Y$ is a \gs{} if and only if the product space $X\x Y$ is. \qed
\ecor

The following observation is made in~\cite[Corollary 24]{Jor07}.

\bprp[Jordan]
Let $X$ be a space.
If the space $\Cp(X)$ is productively \FU{}, then the space $X$ is productively $\gamma$.
\eprp
\bpf
Let $Y$ be a \gs{}. To prove that $X\x Y$ is a \gs{}, we may assume that the spaces $X$ and $Y$ are disjoint.
By the Gerlits--Nagy Theorem, the space $\Cp(Y)$ is \FU{}.
Thus, the space $\Cp(X\sqcup Y)=\Cp(X)\x\Cp(Y)$ is \FU{}.
Applying the Gerlits--Nagy Theorem again, we have that $X\sqcup Y$ is a \gs{}.
Apply Corollary~\ref{cor:unprod}.
\epf

Some of the major results concerning the property that $\Cp(X)$ is productively \FU{}
are collected in the following theorem.

\bthm[Jordan]\label{thm:j}
\mbox{}
\be
\item Assuming \CH{}, there is an uncountable set $X\sub\R$ such that $\Cp(X)$ is productively \FU{}
~\cite[Theorem 33]{Jor07}.
\item There is no uncountable set $X\sub\R$, of cardinality smaller than $\fb$,
such that $\Cp(X)$ is productively \FU{}~\cite[Theorem 34]{Jor07}.
\item The minimal cardinality of a set $X\sub\R$ such that $\Cp(X)$ is not productively \FU{}
is $\aleph_1$~\cite[Corollary 35]{Jor07}.
\item Every uncountable set $X\sub\R$ has a co-countable subset $Y$ such that
$\Cp(Y)$ is not productively \FU{}~\cite[Theorem 1]{Jor08}.
\item If $\Cp(X)$ is productively \FU{}, then so is $\Cp(A)$ for every F$_\sigma$ subset
$A$ of $X$~\cite[Proof of Theorem 1]{Jor08}.
\ee
\ethm

Items (4) and (5) of Jordan's Theorem~\ref{thm:j} solved Problems 1 and 4 of Jordan's
earlier paper~\cite{Jor07}.
The following problem---Problem 3 of~\cite{Jor07}---remains open.

\bprb[{Jordan}]
Is the existence of uncountable set $X\sub\R$ with $\Cp(X)$ productively \FU{}
compatible with Martin's Axiom and the negation of \CH{}?
\eprb

Problem 2 of Jordan~\cite{Jor07} asks whether \CH{} is necessary in item (1).
We solve this problem.
To this end, we use the following characterization of Jordan
\cite[Corollary 23]{Jor07}.
For families of sets $\rmA$ and $\rmB$, let
$$\rmA\land\rmB = \set{\cB\cap \cA}{\cB\in\rmB, \cA\in\rmA}.$$
A family of sets is \emph{centered} if every intersection of finitely many elements from this family is infinite.
A \emph{pseudointersection} of a family $\cF$ of sets is an infinite set $A$ such that $A\as B$ for each
element $B\in\cF$.

\bthm[Jordan]\label{thm:jchar}
Let $X$ be a space, and $\cO$ be the family of all open subsets of $X$. The following two assertions are equivalent:
\be
\item The space $\Cp(X)$ is productively \FU{}.
\item For each family $\rmA\sub\Om(X)$ that is closed under finite intersections, the first property below implies the second:
\bi
\item[(P1)] For every countable family $\rmB\sub P(\cO)$ with $\rmB\land\rmA$ centered,
the family $\rmB\land\rmA$ has a pseudointersection.
\item[(P2)] The family $\rmA$ has a pseudointersection $\cU$ such that $\cU\in\Ga(X)$.
\ei
\ee
\ethm

\blem\label{lem:P1imp}
Let $X$ be a space and $\rmA\sub\Om(X)$ be closed under finite intersections and such that (P1) holds.
Then:
\be
\itm For each countable set $C\sub X$ such that $C$ is not contained in any element of any member of $\rmA$, the family $\rmA$ has a pseudointersection $\cU$ such that $\cU\in\Ga(C)$.
\itm For every sequence $\seq{\cU_n\in P(\cO)}{n\in\N}$ with $\{\cU_n\}\land\rmA$ centered for each $n$,
there is a selection of finite sets $\seq{\cF_n\sub\cU_n}{n\in\N}$
such that the family $\Un_n\cF_n$ is a
pseudointersection of $\rmA$.
\ee
\elem
\bpf
(1) For each finite $F\sub C$, we have that
$$[F]:=\set{U\sub X}{U\mbox{ is open and }F\sub U}\in\Om(X).$$
Let
$$\rmB=\set{[F]}{F\in[C]^{<\oo}}.$$
As $\rmB\land\rmA$ is centered, it has a pseudointersection $\cU$.
In particular, the family $\cU$ is a pseudointersection of $\rmB$, and thus $\cU\in\Ga(C)$.

(2) For each $n$, let $\cV_n=\Un_{m\ge n}\cU_m$. Let $\rmB=\set{\cV_n}{n\in\N}$.
By (P1), the set $\rmB\land\rmA$ has a pseudointersection $\cU$. Represent $\cU=\Un_n\cF_n$ such that
$\cF_n$ is a finite subset of $\cU_n$ for all $n$.
\epf

The following theorem is the main theorem of this section.
Identify $\PN$ with the Cantor space $\Cantor$, via characteristic functions.
The space $\PN$ is homeomorphic to the Cantor set, and can be viewed as a subset of $\R$.
Naturally, the space $\PN$ is the union of $\roth$ and $\Fin$,
the family of infinite subsets of $\N$ and the family of finite subsets of $\N$,
respectively.
We identify elements $x\in\roth$ with increasing elements of $\NN$
by letting $x(n)$ be the $n$th element in the increasing enumeration of $x$.
A subset of $\roth$ is \emph{unbounded} if it is unbounded (with respect to $\le^*$) when
viewed a subset of $\NN$.
An enumerated set $T=\set{x_\alpha}{\alpha<\kappa}$ is a \emph{tower} if
the sequence $\seq{x_\alpha}{\alpha<\kappa}$ is decreasing with respect to $\as$.
Unbounded towers of cardinality $\aleph_1$ exist if and only if $\fb=\aleph_1$
(cf.~\cite[Lemma 3.3]{LinSAdd}).

\bthm\label{thm:gprod}
For each unbounded tower $T=\set{x_\alpha}{\alpha<\aleph_1}$, the space $\Cp(T\cup\Fin)$ is productively \FU{}.
In particular, the space $T\cup\Fin$ is productively $\gamma$.
\ethm
\bpf
Let $X=T\cup\Fin$.
For each $\alpha<\aleph_1$, let $X_\alpha=\sset{x_\beta}{\beta<\alpha}\cup\Fin$.
We may assume that there is $\alpha_0<\aleph_1$ such that $X_{\alpha_0}$ is not contained in any member of any of the considered covers. Indeed,
let $\sset{V_n}{n\in\N}$ be the set of all finite unions of basic open sets.
We may restrict attention to open covers contained in $\sset{V_n}{n\in\N}$.
For each $n$, using that $X$ is not a subset of $V_n$,
let $\beta_n<\aleph_1$ be such that $X_{\beta_n}\not\sub V_n$.
Take $\alpha_0=\sup_n\beta_n$.
Let $\rmA\sub\Om(X)$ be closed under finite intersections and such that (P1) holds.

By (P1) and Lemma~\ref{lem:P1imp}(1), there is a pseudointersection $\cU$ of $\rmA$
such that $\cU\in\Ga(X_{\alpha_0})$.
By~\cite[Lemma 1.2]{GM},
there are $m^0_0<m^0_1<\dots$ and distinct elements $U^0_0,U^0_1,\dots\in\cU$
(so that $\sset{U^0_n}{n\in\N}\in\Ga(X_{\alpha_0})$) such that, for each
$x\in \PN$ and each $n$ with $x\cap (m^0_n,m^0_{n+1})=\emptyset$, we have that $x\in U^0_n$.
Note that $\sset{U^0_n}{n\in\N}$ is a pseudointersection of $\rmA$.
Let $I_0=\N$.

As $\alpha_0<\aleph_1$, the set $\sset{x_\alpha}{\alpha_0<\alpha<\aleph_1}$ is unbounded.
Thus (e.g.,~\cite[Lemma 3.1]{LinSAdd}), there is $\alpha_1>\alpha_0$ such that
the set $I_1:=\sset{n}{x_{\alpha_1}\cap (m^0_n,m^0_{n+1})=\emptyset}$ is infinite.

By (P1) and Lemma~\ref{lem:P1imp}(1), there is a pseudointersection $\cU$ of $\rmA$
such that $\cU\in\Ga(X_{\alpha_1})$.
By~\cite[Lemma 1.2]{GM}, there are $1<m^1_0<m^1_1<\dots$ and distinct elements $U^1_0,U^1_1,\dots\in\cU$
(so that $\sset{U^1_n}{n\in\N}\in\Ga(X_{\alpha_1})$) such that, for each
$x\in \PN$ and each $n$ with $x\cap (m^1_n,m^1_{n+1})=\emptyset$, we have that $x\in U^1_n$.
Here too, the set $\sset{U^1_n}{n\in\N}$ is a pseudointersection of $\rmA$.

Continue in the same manner to define, for each $k>0$, elements with the following properties:
\be
\itm $\alpha_k>\alpha_{k-1}$;
\itm $I_k:=\sset{n}{x_{\alpha_k}\cap (m^{k-1}_n,m^{k-1}_{n+1})=\emptyset}$ is infinite;
\itm $k<m^k_0<m^k_1<\dots$;
\itm $\set{U^k_n}{n\in\N}\in\Ga(X_{\alpha_k})$, and is a bijectively enumerated
pseudointersection of $\rmA$;
\itm For each $x\in \PN$ and each $n$ with $x\cap (m^k_n,m^k_{n+1})=\emptyset$, we have that $x\in U^k_n$.
\ee
Let $\alpha=\sup_k\alpha_k$. Then $\alpha<\aleph_1$, the set $X_\alpha$ is countable, and $X_{\alpha_k}\sub X_{\alpha_{k+1}}$ for all $k$.
Thus, there are for each $k$ a finite set $F_k\sub X_{\alpha_k}$
such that $F_k\sub F_{k+1}$ for all $k$, and $X_\alpha=\Un_kF_k$.
For each $k$, by removing finitely many elements from the set $I_k$, we may assume that
$F_k\sub U^k_n$ for all $n\in I_k$.

Fix $k\in\N$.
By removing finitely many more elements from each set $I_{k+1}$, we may assume that
$x_\alpha\sm[0,m^k_n)\sub x_{\alpha_{k+1}}$ for all $n\in I_{k+1}$.
As $x_{\alpha_{k+1}}\cap (m^k_n,m^k_{n+1})$ is empty for $n\in I_{k+1}$,
we have that
$$x_\alpha\cap (m^k_n,m^k_{n+1})=\emptyset$$
for all $n\in I_{k+1}$.

For each $k$, let $\cU_k=\set{U^k_n}{n\in I_{k+1}}$. By thinning out the sets $I_k$,
we may assume that the families $\cU_k$ are pairwise disjoint.
By Lemma~\ref{lem:P1imp}(2), there are finite
sets $\cF_k\sub\cU_k$ for $k\in\N$ such that $\cU:=\Un_k\cF_k$ is a pseudointersection of $\rmA$.
It remains to show that $\cU\in\Ga(X)$.
Let $x\in X_\alpha$. Let $N$ be such that $x\in F_N$. Then, for each $k\ge N$
and each $U^k_n\in\cF_k$, we have that
$$x\in F_N\sub F_k\sub U^k_n.$$
This shows that $\cU\in\Ga(X_\alpha)$.

It remains to consider the elements $x_\beta$ for $\beta\ge\alpha$.
Let $\beta\ge\alpha$. Then $x_\beta\as x_\alpha$.
Let $k$ be such that $x_\beta\sm [0,k)\sub x_\alpha$.
For each element $U^k_n\in\cF_k$, we have that $n\in I_{k+1}$ and $m^k_n>k$.
Thus,
$$x_\beta\cap (m^k_n,m^k_{n+1})\sub
x_\alpha\cap (m^k_n,m^k_{n+1})=\emptyset,$$
and therefore $x_\beta\in U^k_n$.
\epf

Our proof method cannot produce sets of cardinality greater than $\aleph_1$, since the countability
of the initial sets $X_\alpha$ (for $\alpha<\aleph_1$) is used in an essential manner.

\bcor
The following assertions are equivalent:
\be
\item $\fb=\aleph_1$.
\item There is a set $X\sub\R$, of cardinality $\aleph_1$,
such that $\Cp(X)$ is productively \FU{}.
\ee
\ecor
\bpf
If $\fb=\aleph_1$, then there is an unbounded tower of cardinality $\aleph_1$, and
Theorem~\ref{thm:gprod} applies.
The remaining implication follows from Jordan's Theorem~\ref{thm:j}(2).
\epf

The partial orders $\le^*$ and $\as$, and their inverses, all have the property
mentioned in the following result, that rules out
the possibility of our method to produce examples of cardinality greater than $\aleph_1$.
This is in contrast to~\cite[Theorem 3.6]{LinSAdd}, which implies that \gs{}s $X\sub\R$
of cardinality $\fp$ exist whenever $\fp=\fb$.

\bprp
Assume that $\fb>\aleph_1$. Let $\preceq$ be a partial order on $\roth$ such that,
for each $a\in\roth$, the set $\sset{b\in\roth}{b\preceq a}$ is F$_\sigma$ in $\roth$.
Let $T=\set{x_\alpha}{\alpha<\kappa}$ be strictly $\preceq$-increasing with $\alpha$.
Then the space $\Cp(T\cup\Fin)$ is \emph{not} productively \FU{}.
\eprp
\bpf
Assume that $\Cp(T\cup\Fin)$ is productively \FU{}.

If $\kappa=\aleph_1$, then $\kappa<\fb$ and Jordan's Theorem~\ref{thm:j}(2) applies.
Thus, assume that $\kappa>\aleph_1$.
Let $A=\set{x_\alpha}{\alpha\le\aleph_1}$. As
$$A = T\cap\set{x\in\roth}{x\preceq x_{\aleph_1}},$$
the set $A$ is F$_\sigma$ in $T$.
As $|A\cup\Fin|=\aleph_1<\fb$, the set $A\cup\Fin$ is a $\sigma$-set, that is,
all G$_\delta$ subsets of this set are relatively F$_\sigma$. In particular, the set $A$ is F$_\sigma$ in
$A\cup\Fin$.
Let $F_1$ and $F_2$ be F$_\sigma$ subsets of $\PN$ such that $F_1\cap T=A$ and
$F_2\cap (A\cup Q)=A$. Then
$$F_1\cap F_2\cap (T\cup Q) = (F_1\cap F_2\cap T)\cup (F_1\cap F_2\cap Q) = A\cup \emptyset=A.$$
It follows that $A$ is F$_\sigma$ in $T\cup Q$.
By Jordan's Theorem~\ref{thm:j}(5),
the space $\Cp(A)$ is productively \FU{}, and has cardinality $\aleph_1$, in contradiction to
Jordan's Theorem~\ref{thm:j}(2).
\epf

\bprb\label{p1}
Is the assumption $\fb=\aleph_1$ necessary for the existence of
uncountable sets $X\sub\R$ such that $\Cp(X)$ is productively \FU{}?
\eprb

By Jordan's Theorem~\ref{thm:j}(2),
if the answer to Problem~\ref{p1} is ``No'', then the answer to the following problem
is ``Yes.''

\bprb
Are there, consistently, sets $X\sub\R$ of cardinality greater than $\aleph_1$
such that $\Cp(X)$ is productively \FU{}?
\eprb

\bprb
Are there, consistently, sets $X\sub\R$ such that $X$ is productively $\gamma$
but $\Cp(X)$ is not productively \FU{}?
\eprb


\section{A product of \gs{}s need not have Menger's property}\label{sec:nonprod}

Rothberger's property $\sone(\Op,\Op)$ implies Borel's closely related property of \emph{strong measure zero}. Weiss~\cite{WeissPhD} and, independently, Scheepers~\cite{smzpow} proved that
every metric space satisfying $\ufin(\Op,\Ga)$ and $\sone(\Op,\Op)$ is productively strong measure zero.

\bprb[Scheepers~\cite{smzpow}]\label{SWP}
Assume that $X\sub\R$ satisfies $\ufin(\Op,\Ga)$ and $\sone(\Op,\Op)$.
Must $X$ be productively $\sone(\Op,\Op)$?
\eprb

In~\cite{BabSch}, Babinkostova and Scheepers conjecture that a very strong negative answer to
the Scheepers Problem holds, namely, that assuming \CH{},
there are \gs{}s $X,Y\sub\R$ such that the product space $X\x Y$ does not satisfy $\sfin(\Op,\Op)$.
By Theorem~\ref{thm:gprod}, the unbounded tower method from~\cite{GM, LinSAdd, ProdLind}
cannot be used to establish this conjecture.
Here, we use the Aronszajn tree method of Todor\v{c}evi\'{c}~\cite{GM, Tod95, Brendle96, MilNonGamma}
to prove the Babinkostova--Scheepers Conjecture.

\bthm[CH]\label{thm:nonprod}
There are sets $X,Y\sub\R$ satisfying $\smallbinom{\BO}{\BG}$ such that the product space
$X\x Y$ does not satisfy $\sfin(\Op,\Op)$.
\ethm

In the proof of Theorem~\ref{thm:nonprod}, we work in $\Cantor$ instead of $\R$.
We construct an Aronszajn tree of perfect sets determined by Silver forcing~\cite{grig}.

\bdfn
The partially ordered set $\poset$ is the set of conditions $p$ such that
there is a co-infinite set $D\sub\N$ with $p\colon D\to\two$.
For $p\in\poset$,
$$[p]:=\set{x\in \Cantor}{p\sub x}.$$
A condition $p\in\poset$ is \emph{stronger} than a condition $q\in\poset$, denoted $p\le q$, if
$p\supseteq q$ or, equivalently, if $[p]\sub [q]$.
For $n\in\N$, the relation $p\le_n q$ holds if
$p\le q$ and the first $n$ elements of ${D_p}\comp$ are the same as the first $n$ elements of ${D_q}\comp$.
\edfn

The following important lemma is folklore.

\blem[Fusion Lemma]
Let $\seq{p_n}{n\in\N}$ be a sequence in $\poset$ such that $p_{n+1}\le_n p_n$ for all $n$.
Then the fusion $q=\Un_n p_n$ is in $\poset$, and $q\le_n p_n$ for all $n$.
\elem

\bpf[Proof of Theorem~\ref{thm:nonprod}]
Define the following countable dense subsets of $[p]$:
\begin{eqnarray*}
Q^0(p) & = & \set{x\in [p]}{\forall^\infty n\in{D_p}\comp,\ x(n)=0};\\
Q^1(p) & = & \set{x\in [p]}{\forall^\infty n\in{D_p}\comp,\ x(n)=1}.
\end{eqnarray*}
Define $q\le_n^* p$ if and only if $q\le_n p$ and $q$ is identically zero on $D_q\sm D_p$.

\blem\label{lemone}
Let $\cU\in\Om(Q^0(p))$. For each $n$, there are $U\in\cU$ and $q\le_n^* p$ such that $[q]\sub U$.
\elem
\bpf
Let $F$ be the set consisting of the first $n$ elements of ${D_p}\comp$.
For each $s\in \two^F$, let $x_s\in Q^0(p)$ be such that $x_s\restriction F=s$ and
 $x_s(k)=0$ for every $k\in {D_p}\comp\sm F$.
Take $U\in\cU$ with $\sset{x_s}{s\in \two^F}\sub U$.
Since $U$ is open there is $N\in\N$ with
$[x_s\restriction N]\sub U$ for all $s\in \two^F$.
Define $q\le_n^* p$ by
$$q=p\cup\set{\la k,0\ra}{k<N \mbox{ and } k\in ({D_p}\comp\sm F)}.\qedhere$$
\epf

\blem
Let $p_n\in\poset$, $k_n\in\N$ for $n<N$, and $\cU\in\Om(\Un_{n<N}Q^0(p_n))$.
Then there are $U\in\cU$ and
$\seq{q_n\le_{k_n}^* p_n}{n<N}$ such that
$$\Un_{n<N} [q_n]\sub U.$$
\elem
\bpf
Let $F_n$ be the set consisting of the first $k_n$ elements of ${D_{p_n}}\comp$.
For $s\in \two^{F_n}$, define $x_s^n\in Q^0(p_n)$ as in the
proof of Lemma~\ref{lemone}.
Let $H\sub \Un_{n<N}Q^0(p_n)$ be a finite set containing
all such $x_s^n$. Choose $U\in\cU$ with $H\sub U$ and determine the $q_n$
for $n<N$ as in Lemma~\ref{lemone}.
\epf

\brem
If $q\le_k^* p$ then $Q^0(q)\sub Q^0(p)$ and
hence any $\Om(Q^0(p))\sub\Om(Q^0(q))$.
In these two lemmata, the $q$ we obtain are also equal mod finite
to the $p$, which also implies this.
\erem

\blem\label{gamcover}
Let $\seq{(p_n,k_n)}{n\in\N}$ be a sequence in $\poset\x\N$ and
$\seq{\cU_n}{n\in \N}$ be a sequence in $\Om(Q)$, where $Q=\Un_{n\in\N}Q^0(p_n)$.
Then there are sequences $\seq{U_m\in\cU_m}{m\in \N}$ and
$\seq{q_n\le_{k_n}p_n}{n\in\N}$ such that
$$(\forall n, \forall m\ge n)\quad [q_n]\sub U_m.$$
\elem
\bpf
Construct
$\seq{q_n^m}{n,m\in \N}$ and $\seq{U_m\in\cU_m}{m\in \N}$ by induction on $m$.
Set $q^1_n=p_n$ for all $n$.
Given $\seq{q_n^m}{n\in \N}$
and $\seq{U_n}{n<m}$, construct $q_n^{m+1}$
and $U_m\in\cU_m$ so that
\be
\item $q_n^{m+1}=p_n$ for $n\ge m+1$,
\item $q_n^{m+1}\le_{k_n+m}^* q_n^{m}$ for $n\le m$, and
\item $[q_n^{m+1}]\sub U_m$ for $n\le m$.
\ee
Let $q_n=\Un_{m>n}q_n^m$ be the fusion. We
have that $q_n\le_{k_n}q_n^n=p_n$ and $[q_n]\sub U_m$ whenever
$m\ge n$.
\epf

\brem \label{78}
\mbox{}
\be
\item The analogue of this Lemma for $Q^1$ is also true.
\item The proof of the lemma above only uses
the fact that $[p_n]\cap U$ is open in $[p_n]$ for all $n$
and $U$ appearing in some $\cU_m$.
\ee
\erem

\blem\label{borel}
Let $p\in\poset$, $n\in\N$, and $B\sub \Cantor$ be a Borel set.
Then there exists
$q\le_n p$ such that $[q]\cap B$ is clopen
in $[q]$.
\elem
\bpf
Let $F$ be the set consisting of the first $n$ elements of ${D_p}\comp$
and let $\phi\colon \N\to ({D_p}\comp\sm F)$ be a bijection.
For $I\sub\N$ let $\psi_I\colon ({D_p}\comp\sm F)\to\two$ be
the restriction of the characteristic function of $\phi(I)$.
For each $s\in \two^F$ define
$$C_s=\set{I\in\roth}{(p\cup s\cup \psi_I)\in B}.$$
Since these are Borel sets,
by the Galvin--Prikry Theorem~\cite{gp} there exists
$H\in\roth$ such that for each $s\in \two^F$ either
$[H]^\infty\sub C_s$ or $[H]^\infty\cap C_s=\emptyset$.
Let $H_1\sub H$ be infinite such that $H\sm H_1$ is also infinite.
Let
$$q=p\cup (\phi(H\comp)\x\{0\})\cup (\phi(H_1)\x\{1\}).$$
Note that ${D_q}\comp=F\cup \phi(H\sm H_1)$.
We claim that given any $x,y\in [q]$,
if $x\restriction F= y \restriction F = s$, then $x\in B$ if and only if $y\in B$.
Letting $H_x=\phi^{-1}(x^{-1}(1))$, we have that
$H_1\sub H_x\sub H$ and so $H_x$ is an infinite subset of $H$.
Similarly for $H_y$. By the choice of $H$ we have that
$H_x\in C_s$ if and only if $H_y\in C_s$, and the claim follows.
\epf

\blem\label{disj}
Let $\seq{(p_n,k_n)}{n\in\N}$ be a sequence in $\poset\x\N$.
Then there is a sequence $\seq{q_n\le_{k_n}p_n}{n\in\N}$ such that for
$n\neq m$,
$q_n$ and $q_m$ are strongly disjoint, i.e., there
are infinitely many $k\in (D_{q_n}\cap D_{q_m})$ with
$q_n(k)\neq q_m(k)$.
\elem
\bpf
Given $p_1,p_2$ and $n$ it is easy to find $q_1\le_n p_1$ and
$q_2\le_n p_2$ which are strongly disjoint. A fusion
argument produces a sequence $\seq{q_n}{n\in\N}$ where all pairs
have been considered and made strongly disjoint.
\epf

We construct an Aronszajn tree of Silver conditions.
Let $B^\beta$ for $\beta<\aleph_1$ list all Borel sets.
Let $\cB_\alpha=\seq{\cB_\alpha^n}{n\in\N}$ for $\alpha<\aleph_1$ be all countable sequences of
countable families of Borel sets.
We may assume that each
element of $\Un_n \cB_\alpha^n$ is equal to
$B^\beta$ for some $\beta<\alpha$.
We may also assume that each such sequence occurs as an element $\cB_\alpha$ for
both $\alpha$ even and $\alpha$ odd.

We construct a tree $T\sub \N^{<\aleph_1}$ and
$\seq{p_s\in\poset}{s\in T}$ with the following properties:
\be
\item $T\sub\N^{<\aleph_1}$ is a subtree, i.e.,
$s\sub t\in T$ implies $s\in T$.
\item $T_\alpha=T\cap\N^\alpha$ is countable for each $\alpha<\aleph_1$.
\item $s\sub t\in T$ implies $p_t\le p_s$.
\item If $s,t\in T$ are incomparable, then $p_s$ and
$p_t$ are strongly disjoint (as in Lemma~\ref{disj}).
\item For any $\alpha<\beta<\aleph_1$ and any $s\in T_\alpha$ and
$n\in\N$ there is $t\in T_\beta$ with $p_t\le_n p_s$.
\item For any $\beta<\alpha$ and $s\in T_\alpha$,
$[p_s]\cap B^\beta$ is clopen in $[p_s]$.
\item Define
\begin{eqnarray*}
Q^0_\alpha & = & \Un\set{Q^0(p_t)}{t\in T_{\le\alpha}};\\
Q^1_\alpha & = & \Un\set{Q^1(p_t)}{t\in T_{\le\alpha}}.
\end{eqnarray*}
\be
\item
For an even ordinal $\alpha$, if $\cB_\alpha=\seq{\cB_n^\alpha}{n\in\N}$ is
a sequence in $\Om(Q^0_\alpha)$ then
there is a family
$$\seq{U_n\in\cU_n^\alpha}{n\in\N}\in\Ga(Q^0_\alpha\cup\Un\set{[p_s]}{s\in T_{\alpha+1}}).$$
\item For $\alpha$ odd, the analogous statement
is true with $Q^1_\alpha$ in place of $Q^0_\alpha$.
\ee
\item Let $D=\sset{ {D_{p_s}}\comp}{s\in T}\sub\roth$.
Then $D$ is dominating.
\ee
To construct $T_\lambda$ and $p_s$ for $s\in T_\lambda$ where
$\lambda$ is a countable limit ordinal, proceed as follows.
For any $s\in T_{<\lambda}$ and $N\in\N$
choose a strictly increasing sequence $\seq{\lambda_n}{n\in\N}$ cofinal in $\lambda$ with
$s\in T_{\lambda_1}$. Let $t_1=t^{s,N}_1$ be equal to $s$.
By the inductive hypothesis we can find
$t_n=t^{s,N}_n\in T_{\lambda_n}$ with $p_{t_{n+1}}\le_{N+n} p_{t_n}$ for all $n$.
Set $t^{s,N}=\Un_{n}t^{s,N}_n$ and $T_\lambda=\sset{t^{s,N}}{s\in T_{<\lambda}, N\in\N}$.
For every $t=t^{s,N}\in T_\lambda$, let $p_t$ be the fusion of
the sequence $\seq{p_{t^{s,N}_n}}{n\in\N}$, i.e.,
$p_t=\Un_n p_{t^{s,N}_n}$.

At successor stages for $\alpha$ even, check to see if $\cB_\alpha$ is a
sequence in $\Om(Q^0_\alpha)$. If it is not,
we need never worry about
it since the set we are building will contain $Q^0_\alpha$.
If it is, let $\sset{x_n}{n\in\N}=Q^0_\alpha$ and let
$$\cB_n=\set{B\in\cB_n^\alpha}{\set{x_i}{i<n}\sub B}.$$
Let $\seq{p_n,k_n}{n\in\N}$ list all elements of
$$\set{p_s}{s\in T_\alpha}\x\N$$
with infinite repetitions.
Combining the fact that only $B^\beta$'s for $\beta<\alpha$ may occur in some $\cB_n^\alpha$,
Lemma~\ref{gamcover} (see also Remark~~\ref{78}), and
 Lemma~\ref{disj},
 we can find sequences
$\seq{q_n\le_{k_n} p_n}{n\in\N}$ and $\seq{B_m\in\cB_m}{m\in\N}$
such that $[q_n]\sub B_m$ for all $n<m$ and $q_{n_1}, q_{n_2}$ are strongly disjoint
for all distinct $n_1,n_2\in \N$.
As a result, for every $s\in T_\alpha$ and $k\in\N$ there is
some $q_{s,k}\le_{k} s$ such that $[q_{s,k}]\sub B_m$ for all but finitely many $m$.
By Lemma~\ref{borel}, for such $s$ and $k$ there is $p\le_k q_{s,k}$
such that $[p]\cap B^\alpha$ is clopen in $[p]$. We denote this $p$ by $p_{s\concat{k}}$.

This concludes our inductive construction, which ensures Conditions (1)--(7).
Obtaining Condition (8) is easy to satisfy.
Set
\begin{eqnarray*}
X & = & \Un_{s\in T} Q^0(p_s);\\
Y & = & \Un_{s\in T} Q^1(p_s).
\end{eqnarray*}
By Condition (7), the sets $X$ and $Y$ satisfy
$\smallbinom{\BO}{\BG}$.
For all $x\in X$ and $y\in Y$,
there are infinitely many $n$ with $x(n)\neq y(n)$.
Indeed, if $x\in Q^0(p_s)$ and $y\in Q^1(p_t)$,
and $s$ and $t$ are incomparable, then $p_s$ and
$p_t$ are strongly disjoint. On the other hand, if $s$ and $t$
are comparable, for example, if $s\sub t$, then
since $p_t\le p_s$, we have that ${D_{p_t}}\comp\sub {D_{p_s}}\comp$.
Thus, for all but finitely many $n\in {D_{p_t}}\comp$, we have
that $y(n)=1$ and $x(n)=0$.

Condition (8) provides a continuous map from $X\x Y$ onto
a dominating set $D\sub\NN$. Namely, if
$x_0\in Q^0(p_s)$ is identically zero on ${D_{p_s}}\comp$
and $x_1\in Q^1(p_s)$ is identically one on ${D_{p_s}}\comp$,
then ${D_{p_s}}\comp=\sset{n}{x_0(n)\neq x_1(n)}$.
Thus, the continuous map
$\Phi\colon X\x Y\to \NN$
defined by $\Phi(x,y)=\sset{n}{x(n)\neq y(n)}$, is as required.
\epf


\section{Applications}\label{sec:apps}

The conjunction of Hurewicz's property $\ufin(\Op,\Ga)$ and Rothberger's property
$\sone(\Op,\Op)$, shown in~\cite[Theorems 14 and 19]{NSW} to be equivalent to
the Gerlits--Nagy property
(*), is of growing importance in the area of selection principles~\cite{SchGerlits}.
In an unpublished manuscript~\cite{WeissProds}, Weiss proposed a plan to prove that
the Gerlits--Nagy property (*) is preserved by finite products.
By Lemma~\ref{lem:un2prod}, this problem is equivalent to the following one.

\bprb[Weiss]\label{prb:weiss}
Is the conjunction of $\ufin(\Op,\Ga)$ and $\sone(\Op,\Op)$ preserved by finite powers?
\eprb

A negative solution of Weiss's Problem was proposed in~\cite{SchTall10}, and later withdrawn~\cite{STErrata}.
A set $S\sub\R$ is \emph{Sierpi\'nski} if the set $S$ is uncountable, and its intersection with every
Lebesgue measure zero set is countable.
The solution proposed in~\cite{SchTall10} was based on the assumption that if $S\sub\R$ is a Sierp\'nski set,
then $S$ continues to satisfy the Hurewicz property $\ufin(\Op,\Ga)$
in extensions of the universe by Cohen forcing~\cite[Theorem 40]{SchTall10}.
It turns out that this assumption is not provable
(Theorem~\ref{thm:starpres} below).\footnote{
The gap in the proof of Theorem 40 in~\cite{SchTall10} may be the following one.
It seems that, in item 6) on page 30, the definition of $\name{V}^n_j$ should be
$\name{V}^n_{j-1}\cap\bigl(\bigcap_{i\le \ell^n_j} \name{V}^n_{m^n_i+\cdots+m^n_{j-1}+1, x^{n, j}_i}\bigr)$,
not
$\name{V}^n_{j-1}\cap\bigl(\bigcap_{i\le \ell^n_j} \name{V}^n_{j, x^{n, j}_i}\bigr)$.
Given that, the claim
``By 3), 5), 6) and 8) above, the set $F_k$ is disjoint from
$\Un_{n\ge k}C_n$'' at the end of page 30 is unclear.
Indeed, to make it true,
one should have in $V[G]$ that $\overline{V^n_t}\spst C_n$.
By the definition of $V^n_j$, this would require
that, in $V[G]$,
$C_n\sub \overline{V^n_{m^n_i+\cdots+m^n_{j-1}+1, x^{n, j}_i}}$
for all $i<\ell^n_j$.
For each individual $i<\ell^n_j$, every element
$p$ of $F_{m^n_i+\cdots+m^n_{j-1}+1}(x^{n,j}_i, \name{C}_n)$
indeed forces that $\name{C}_n\sub \name{V}_{p,x^{n,j}_i}(\name{C}_n)$.
However, the elements of
$F_{m^n_i+\cdots+m^n_{j-1}+1}(x^{n,j}_i, \name{C}_n)$ may be
incompatible.
As, in $V[G]$, we have that
$$V^n_{m^n_i+\cdots+m^n_{j-1}+1, x^{n, j}_i}=\bigcap \smallmedset{ V_{p,x^{n,j}_i}(\name{C}_n)}
{p\in F_{m^n_i+\cdots+m^n_{j-1}+1}(x^{n,j}_i, \name{C}_n)},$$
it is unclear why
$C_n\sub V^n_{m^n_i+\cdots+m^n_{j-1}+1, x^{n, j}_i}$ there.
}

The following Theorem~\ref{thm:manysol} provides an alternative solution to Weiss's Problem,
also in the negative.
In particular, the answer to Problem 6.6 in~\cite{OPiT} is ``No.''
It was, thus far, open whether the Gerlits--Nagy property (*) implies $\sone(\Om,\Om)$.
Theorem~\ref{thm:manysol} solves this problem, in the negative.
It also shows that the answer to Problem 4.1(j) in~\cite{OPiT}, concerning the realization of
a certain setting in the Borel version of the Scheepers Diagram is ``Yes''. This theorem solves
8 out of the 55 problems that remained open in Mildenberger--Shelah--Tsaban~\cite{MShT858}, concerning
potential implications between covering properties (details are provided below).
It also solves, in the negative, all 5 problems in~\cite[Problem 7.6(2)]{OPiT}, concerning the preservation
of certain covering properties under finite powers.

\begin{figure}[!ht]
\renewcommand{\sr}[2]{{\txt{$#1$\\$#2$}}}
{\tiny
\begin{changemargin}{-3cm}{-3cm}
\begin{center}
$\xymatrix@C=7pt@R=6pt{
&
&
& \sr{\ufin(\Op,\Gamma)}{\fb~~ (18)}\ar[r]
& \sr{\ufin(\Op,\Tau)}{\max\{\fb,\fs\}~~ (19)}\ar[rr]
&
& \sr{\ufin(\Op,\Omega)}{\fd~~ (20)}\ar[rrrr]
&
&
&
& \sr{\sfin(\Op,\Op)}{\fd~~ (21)}
\\
&
&
& \sr{\sfin(\Gamma,\Tau)}{{\fb}~~ (12)}\ar[rr]\ar[ur]
&
& \sr{\sfin(\Gamma,\Omega)}{\fd~~ (13)}\ar[ur]
\\
\sr{\sone(\Gamma,\Gamma)}{\fb~~ (0)}\ar[uurrr]\ar[rr]
&
& \sr{\sone(\Gamma,\Tau)}{{\fb}~~ (1)}\ar[ur]\ar[rr]
&
& \sr{\sone(\Gamma,\Omega)}{\fd~~ (2)}\ar[ur]\ar[rr]
&
& \sr{\sone(\Gamma,\Op)}{\fd~~ (3)}\ar[uurrrr]
\\
&
&
& \sr{\sfin(\Tau,\Tau)}{{\min\{\fs,\fb\}}~~ (14)}\ar'[r][rr]\ar'[u][uu]
&
& \sr{\sfin(\Tau,\Omega)}{\fd~~ (15)}\ar'[u][uu]
\\
\sr{\sone(\Tau,\Gamma)}{\fp~~ (4)}\ar[rr]\ar[uu]
&
& \sr{\sone(\Tau,\Tau)}{{\fp}~~ (5)}\ar[uu]\ar[ur]\ar[rr]
&
& \sr{\sone(\Tau,\Omega)}{\od~~ (6)}\ar[uu]\ar[ur]\ar[rr]
&
& \sr{\sone(\Tau,\Op)}{\od~~ (7)}\ar[uu]
\\
&
&
& \sr{\sfin(\Omega,\Tau)}{\fp~~ (16)}\ar'[u][uu]\ar'[r][rr]
&
& \sr{\sfin(\Omega,\Omega)}{\fd~~ (17)}\ar'[u][uu]
\\
\sr{\sone(\Omega,\Gamma)}{\fp~~ (8)}\ar[uu]\ar[rr]
&
& \sr{\sone(\Omega,\Tau)}{\fp~~ (9)}\ar[uu]\ar[ur]\ar[rr]
&
& \sr{\sone(\Omega,\Omega)}{\cov(\cM)~~ (10)}\ar[uu]\ar[ur]\ar[rr]
&
& \sr{\sone(\Op,\Op)}{\cov(\cM)~~ (11)}\ar[uu]
}$
\end{center}
\end{changemargin}
}
\caption{The Extended Scheepers Diagram}\label{fig:ESD}
\end{figure}
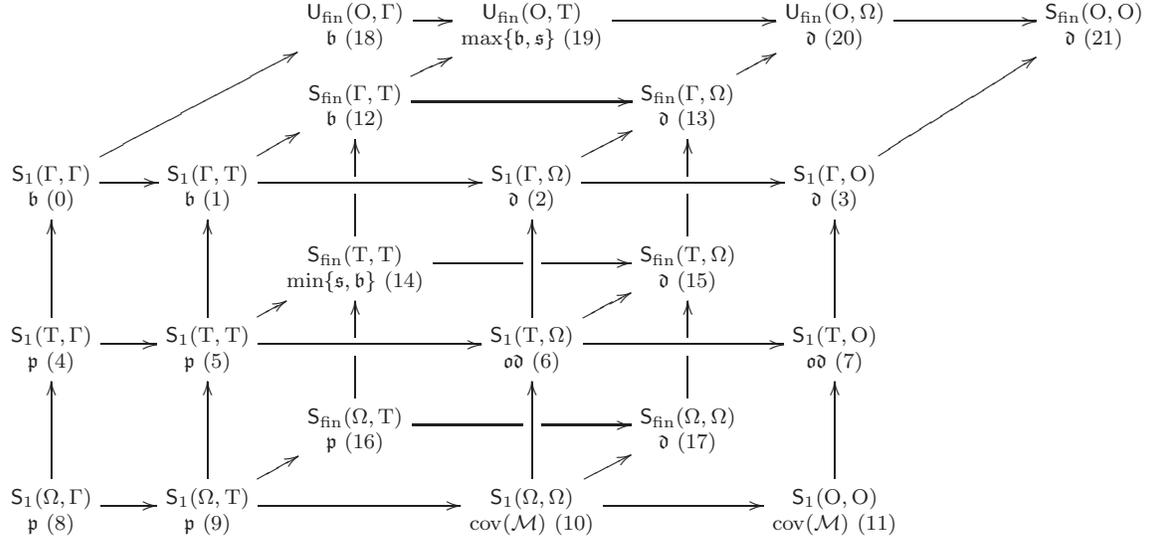

An element $\cU\in\Op(X)$ is in $\Tau(X)$ if every member of $X$ is a member of infinitely many
elements of $\cU$, and, for all $x,y\in X$, either $x\in U$ implies $y\in U$ for all but finitely many $U\in\cU$, or $y\in U$ implies $x\in U$ for all but finitely many $U\in\cU$.
Figure~\ref{fig:ESD} contains all new properties introduced by the
inclusion of $\Tau$ into the framework, together with their critical cardinalities
\cite{tautau, ShTb768, MShT858, MalliarisShelah}, and a serial number to be used below.

\bthm[CH]\label{thm:manysol}
There are sets $X_0,X_1\sub\R$ satisfying $\smallbinom{\BO}{\BG}$,
such that the set $X=X_0\cup X_1$ has the following properties:
\be
\item $X$ satisfies $\sone(\BT,\BG)$ and $\sone(\B,\B)$ (and, in particular, the Gerlits--Nagy property (*));
\item $X$ does not satisfy $\sfin(\Om,\Om)$;
\item The square space $X^2$ does not satisfy $\sfin(\Op,\Op)$.
\ee
\ethm
\bpf
Let $X_0,X_1\sub\R$ be as in Theorem~\ref{thm:nonprod},
i.e.,
both satisfying $\smallbinom{\BO}{\BG}$, and such that the product space
$X_0\x X_1$ does not satisfy $\sfin(\Op,\Op)$.
We may assume, by taking a homeomorphic image, that $X_0\sub (0,1)$ and $X_1\sub (2,3)$.
Let $X=X_0\cup X_1$.

(1) As both properties $\sone(\BT,\BG)$ and $\sone(\rmB,\rmB)$ are preserved by finite unions
(e.g.,~\cite{AddQuad}), $X$ satisfies $\sone(\BT,\BG)$ and $\sone(\rmB,\rmB)$.

(2) This follows from (3), since $\sfin(\Om,\Om)$ is equivalent to being $\sfin(\Op,\Op)$ in all
finite powers~\cite[Theorem 3.9]{coc2}.

(3) The product space $X_0\x X_1$ is closed in $X^2$.
Since Menger's property $\sfin(\Op,\Op)$ is hereditary for closed subsets,
the space $X^2$ does not satisfy $\sfin(\Op,\Op)$.
\epf

\forget
The following diagram indicates exactly which properties in the Scheepers Diagram
are satisfied by the space $X$ of Corollary~\ref{cor:setting}.
$$\xymatrix@R=8pt{
&
& \bullet\ar[r]
& \bullet\ar[rr]
& & \bullet
\\
&
& \bullet\ar[ur]
\\
\bullet\ar[r]\ar[uurr]
& \bullet\ar[rr]\ar[ur]
& & \bullet\ar[uurr]
\\
&
& \circ\ar'[u][uu]
\\
\circ\ar[r]\ar[uu]
& \circ\ar[uu]\ar[rr]\ar[ur]
& & \bullet\ar[uu]
}$$
\forgotten

The set in Theorem~\ref{thm:manysol} realizes the following setting in the Extended Scheepers Diagram.
$$
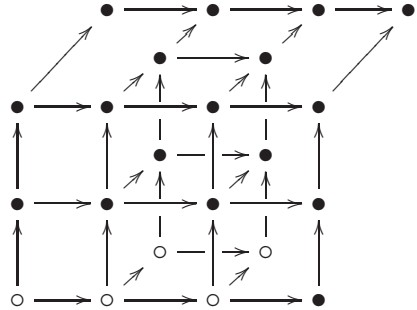
\xymatrix@C=7pt@R=6pt{
&
& \bullet\ar[rr]
&
& \bullet\ar[rr]
&
& \bullet\ar[rr]
&
& \bullet
\\
&
&
& \bullet\ar[rr]\ar[ur]
&
& \bullet\ar[ur]
\\
\bullet\ar[uurr]\ar[rr]
&
& \bullet\ar[ur]\ar[rr]
&
& \bullet\ar[ur]\ar[rr]
&
& \bullet\ar[uurr]
\\
&
&
& \bullet\ar'[r][rr]\ar'[u][uu]
&
& \bullet\ar'[u][uu]
\\
\bullet\ar[rr]\ar[uu]
&
& \bullet\ar[uu]\ar[ur]\ar[rr]
&
& \bullet\ar[uu]\ar[ur]\ar[rr]
&
& \bullet\ar[uu]
\\
&
&
& \circ\ar'[u][uu]\ar'[r][rr]
&
& \circ\ar'[u][uu]
\\
\circ\ar[uu]\ar[rr]
&
& \circ\ar[uu]\ar[ur]\ar[rr]
&
& \circ\ar[uu]\ar[ur]\ar[rr]
&
& \bullet\ar[uu]
}$$
Consider the serial numbers in the Extended Scheepers Diagram. The table below describes all known
implications and nonimplications among the properties, so that entry $(i,j)$ indicates whether property $(i)$ implies
property $(j)$. The framed entries remained open in~\cite{MShT858}. Their solution follows from
Theorem~\ref{thm:manysol}. This gives a complete
understanding of which properties in the Extended Scheepers Diagram imply $\sfin(\Om,\Om)$ and which properties are implied by $\sfin(\Tau,\Om)$.

\begin{table}[!ht]
\begin{changemargin}{-3cm}{-3cm}
\begin{center}
{\tiny
\begin{tabular}{|r||cccccccccccccccccccccc|}
\hline
 & \smb{0} & \smb{1} & \smb{2} & \smb{3} & \smb{4} & \smb{5} & \smb{6} & \smb{7} &
 \smb{8} & \smb{9} & \smb{10} & \smb{11} & \smb{12} & \smb{13} & \smb{14} & \smb{15} &
 \smb{16} & \smb{17} & \smb{18} & \smb{19} & \smb{20} & \smb{21}\cr
\hline\hline

\mb{ 0} &
\yup&\yup&\yup&\yup&\nop&\nop&\nop&\nop&\nop&\nop&\nop&
\nop&\yup&\yup&\nop&\mbq&\nop&\nop&\yup&\yup&\yup&\yup\\
\mb{ 1} &
\mbq&\yup&\yup&\yup&\nop&\nop&\nop&\nop&\nop&\nop&\nop&
\nop&\yup&\yup&\nop&\mbq&\nop&\nop&\mbq&\yup&\yup&\yup\\
\mb{ 2} &
\nop&\nop&\yup&\yup&\nop&\nop&\nop&\nop&\nop&\nop&\nop&
\nop&\nop&\yup&\nop&\mbq&\nop&\nop&\nop&\nop&\yup&\yup\\
\mb{ 3} &
\nop&\nop&\nop&\yup&\nop&\nop&\nop&\nop&\nop&\nop&\nop&
\nop&\nop&\nop&\nop&\nop&\nop&\nop&\nop&\nop&\nop&\yup\\
\mb{ 4} &
\yup&\yup&\yup&\yup&\yup&\yup&\yup&\yup&\nop&\nop&\fbn&
\mbq&\yup&\yup&\yup&\yup&\nop&\fbn&\yup&\yup&\yup&\yup\\
\mb{ 5} &
\mbq&\yup&\yup&\yup&\mbq&\yup&\yup&\yup&\nop&\nop&\fbn&
\mbq&\yup&\yup&\yup&\yup&\nop&\fbn&\mbq&\yup&\yup&\yup\\
\mb{ 6} &
\nop&\nop&\yup&\yup&\nop&\nop&\yup&\yup&\nop&\nop&\fbn&
\mbq&\nop&\yup&\nop&\yup&\nop&\fbn&\nop&\nop&\yup&\yup\\
\mb{ 7} &
\nop&\nop&\nop&\yup&\nop&\nop&\nop&\yup&\nop&\nop&\nop&
\mbq&\nop&\nop&\nop&\nop&\nop&\nop&\nop&\nop&\nop&\yup\\
\mb{ 8} &
\yup&\yup&\yup&\yup&\yup&\yup&\yup&\yup&\yup&\yup&\yup&
\yup&\yup&\yup&\yup&\yup&\yup&\yup&\yup&\yup&\yup&\yup\\
\mb{ 9} &
\mbq&\yup&\yup&\yup&\mbq&\yup&\yup&\yup&\mbq&\yup&\yup&
\yup&\yup&\yup&\yup&\yup&\yup&\yup&\mbq&\yup&\yup&\yup\\
\mb{10} &
\nop&\nop&\yup&\yup&\nop&\nop&\yup&\yup&\nop&\nop&\yup&
\yup&\nop&\yup&\nop&\yup&\nop&\yup&\nop&\nop&\yup&\yup\\
\mb{11} &
\nop&\nop&\nop&\yup&\nop&\nop&\nop&\yup&\nop&\nop&\nop&
\yup&\nop&\nop&\nop&\nop&\nop&\nop&\nop&\nop&\nop&\yup\\
\mb{12} &
\mbq&\mbq&\mbq&\mbq&\nop&\nop&\nop&\nop&\nop&\nop&\nop&
\nop&\yup&\yup&\nop&\mbq&\nop&\nop&\mbq&\yup&\yup&\yup\\
\mb{13} &
\nop&\nop&\nop&\nop&\nop&\nop&\nop&\nop&\nop&\nop&\nop&
\nop&\nop&\yup&\nop&\mbq&\nop&\nop&\nop&\nop&\yup&\yup\\
\mb{14} &
\mbq&\mbq&\mbq&\mbq&\nop&\nop&\nop&\nop&\nop&\nop&\nop&
\nop&\yup&\yup&\yup&\yup&\nop&\fbn&\mbq&\yup&\yup&\yup\\
\mb{15} &
\nop&\nop&\nop&\nop&\nop&\nop&\nop&\nop&\nop&\nop&\nop&
\nop&\nop&\yup&\nop&\yup&\nop&\fbn&\nop&\nop&\yup&\yup\\
\mb{16} &
\mbq&\mbq&\mbq&\mbq&\mbq&\mbq&\mbq&\mbq&\mbq&\mbq&\mbq&
\mbq&\yup&\yup&\yup&\yup&\yup&\yup&\mbq&\yup&\yup&\yup\\
\mb{17} &
\nop&\nop&\nop&\nop&\nop&\nop&\nop&\nop&\nop&\nop&\nop&
\nop&\nop&\yup&\nop&\yup&\nop&\yup&\nop&\nop&\yup&\yup\\
\mb{18} &
\nop&\nop&\nop&\nop&\nop&\nop&\nop&\nop&\nop&\nop&\nop&
\nop&\nop&\mbq&\nop&\mbq&\nop&\nop&\yup&\yup&\yup&\yup\\
\mb{19} &
\nop&\nop&\nop&\nop&\nop&\nop&\nop&\nop&\nop&\nop&\nop&
\nop&\nop&\mbq&\nop&\mbq&\nop&\nop&\nop&\yup&\yup&\yup\\
\mb{20} &
\nop&\nop&\nop&\nop&\nop&\nop&\nop&\nop&\nop&\nop&\nop&
\nop&\nop&\mbq&\nop&\mbq&\nop&\nop&\nop&\nop&\yup&\yup\\
\mb{21} &
\nop&\nop&\nop&\nop&\nop&\nop&\nop&\nop&\nop&\nop&\nop&
\nop&\nop&\nop&\nop&\nop&\nop&\nop&\nop&\nop&\nop&\yup\\

\hline
\end{tabular}
}
\end{center}
\end{changemargin}
\caption{Known implications and nonimplications}\label{imptab}
\end{table}


\section{Preservation under forcing extensions}\label{sec:forcing}

Scheepers proved in~\cite{SchGerlits} that random real forcing preserves
being a \gs{}.
We will show that this is also the case for Cohen's forcing.
We say that a property is \emph{preserved by Cohen forcing} if,
whenever a space $X$ has this property in the ground model, it will have this property in
any extension by Cohen forcing, adding any number of Cohen reals.

\bthm\label{thm:gpres}
The property $\gamma$ is preserved by Cohen forcing.
\ethm
\bpf
Let $M$ be the ground model, and $X$ be a \gs{} in $M$.
Let $G$ be $\poset$-generic over $M$, and $\kappa>0$
be an arbitrary, possibly finite, cardinal.
Let $\poset$ be the poset adding $\kappa$ Cohen reals.
In $M[G]$, let $\cU\in\Om(X)$ be a cover consisting of open sets in $M$.

According to Lemma 3.3 of~\cite{Dow88}, the Lindel\"of property is preserved by
adding uncountably many Cohen reals. The proof of that Lemma also shows that
the Lindel\"of property is preserved by adding countably many Cohen reals.
Thus, in $M[G]$, all finite powers of $X$ are Lindel\"of, and therefore $\cU$ contains a
countable member of $\Om(X)$. Thus, we may assume that $\cU$ is countable, and hence
is determined in an extension by countably many Cohen reals.
As the poset for adding countably many Cohen reals is countable, it is isomorphic to
$\two^{<\alephes}$.
Thus, we may assume that $\poset=\two^{<\alephes}$.
Let $p_0\in\poset$ be a condition forcing the above-mentioned properties of $\cU$.
To simplify our notation, assume that $p_0$ is the trivial condition,
or replace $\poset$ by the conditions stronger than $p_0$.
Work in $M$.

Fix $p\in\poset$. Let
$$\cU_p=\set{U}{\exists q\le p,\ q\forces U\in \name{\cU}}.$$
Then $\cU_p\in\Om(X)$.
As $X$ is a \gs{}, we may, by thinning out $\cU_p$,
assume that $\cU_p\in\Ga(X)$. Thus, by further thinning out, we
may assume that the sets $\cU_p$, for $p\in\poset$, are
pairwise disjoint.
As $X$ satisfies $\sone(\Om,\Ga)$ (the property $\sone(\Ga,\Ga)$ suffices here),
there are elements $U_p\in\cU_p$ for $p\in\poset$
such that $\sset{U_p}{p\in\poset}\in\Ga(X)$.
As the families $\cU_p$ are pairwise disjoint, the sets $U_p$ are distinct for
distinct conditions $p\in\poset$.
For each $p\in\poset$, pick a condition $q_p\le p$ forcing that
$U_p\in \name{\cU}$.

As the set $\sset{q_p}{p\in\poset}$ is dense in $\poset$, its intersection with $G$ is infinite.
Thus, the family $\set{U_p}{q_p\in G}$, which is a subset of $\cU$,
is infinite. As $\sset{U_p}{p\in\poset}\in\Ga(X)$, we have that $\set{U_p}{q_p\in G}\in\Ga(X)$.
\epf

In~\cite[Theorem 37]{SchTall10}, Scheepers and Tall show that the negation of
Hurewicz's property $\ufin(\Op,\Ga)$ is preserved by Cohen forcing.
In~\cite[page 26]{SchTall10}, it is shown that adding a Cohen real destroys the property
that the ground model's Cantor set satisfies $\ufin(\Op,\Ga)$.
Problem 6 in~\cite{SchTall10} asks whether $\ufin(\B,\BG)$---the Hurewicz property for countable
Borel covers---is preserved by Cohen forcing. The following theorem shows, in particular,
that the answer is ``No.''
It is well known that Sierp\'nski sets, which have positive outer measure,
satisfy $\ufin(\B,\BG)$. (A simple proof is given, e.g., in~\cite{MHP}.)
As Rothberger's property $\sone(\Op,\Op)$ implies Lebesgue measure zero,
Sierp\'nski sets cannot satisfy $\sone(\Op,\Op)$.

In the proof of our theorem, we use a technical lemma,
whose proof applies the Rothberger game $\gone(\Op,\Op)$.
This is a game for two players, ONE and TWO, with an inning per each natural number $n$.
In the $n$-th inning, ONE picks a cover $\cU_n\in\Op(X)$, and TWO responds by picking
an element $U_n\in\cU_n$. ONE wins if $\sset{U_n}{n\in\N}$ is not a cover of $X$. Otherwise,
TWO wins. Pawlikowski proved in~\cite[Theorem 1]{Paw} that, for spaces $X$ with points G$_\delta$,
the space $X$ satisfies $\sone(\Op,\Op)$ if and only if ONE does not have a winning strategy in
the game $\gone(\Op,\Op)$.

\bthm\label{thm:starpres}
For $\ufin(\Op,\Ga)$ spaces $X$ with points G$_\delta$, the following assertions are equivalent:
\be
\item $X$ remains $\ufin(\Op,\Ga)$ in every forcing extension by adding Cohen reals.
\item $X$ remains $\ufin(\Op,\Ga)$ in every forcing extension by adding one Cohen real.
\item $X$ satisfies $\sone(\Op,\Op)$.
\ee
\ethm
\bpf
The implication ``$(1)\Impl (2)$'' is trivial. The implication ``$(2)\Impl (1)$''
is proved as in the proof of Theorem~\ref{thm:gpres}, namely, a counter-example
to $\ufin(\Op,\Ga)$ in the extension is determined in an extension by a single Cohen real,
and the negation of $\ufin(\Op,\Ga)$ is preserved by Cohen forcing~\cite[Theorem 37]{SchTall10}.

$(2)\Impl (3)$: Let $M$ be the ground model.
The property $\ufin(\Op,\Ga)$ implies, in particular, that the space $X$ is Lindel\"of in $M$.
Let $\poset=\N^{<\alephes}$, the poset adding one Cohen real $g\in\NN$.

Let $\seq{\cU_n}{n\in\N}\in M$ be a sequence of open covers of $X$.
Since $X$ is Lindel\"of, we may assume that, for each $n$, there is an enumeration
$\cU_n=\set{U^n_m}{m\in\N}$.
Let $G$ be $\poset$-generic over $M$, and $g=\bigcup G\in\NN$ be the corresponding Cohen real.
By genericity, the family $\sset{U^n_{g(n)}}{n\ge k}$ is a cover of $X$ for each $k$.
If the family $\sset{U^n_{g(n)}}{n\in\N}$ has a finite subcover $\sset{U^n_{g(n)}}{n<k}$,
then (since the restriction of $g$ to $\{0,\dots,k-1\}$ is in $M$) this finite subcover is in $M$,
and we are done. Thus, assume that this is not the case.

By (2), there is a function $f\in\NN\cap M[G]$ such that
$$\medset{\bigcup_{k\le n<f(k)}U^n_{g(n)}}{k\in \N}\in\Gamma(X).$$
Work in the ground model.
For $p\in\poset$ and $K\in\N$, let
$$X(p,K)=\medset{x\in X}{p\forces\forall k\ge K,\ x\in \bigcup_{k\le n<\name{f}(k)}U^n_{\name{g}(n)}}.$$
Then $X=\bigcup_{(p,K)\in\poset\x\N}X(p,K)$, a countable union.
We may assume that, for each $k$, $\cU_{k+1}$ is a refinement of $\cU_k$.

\bclm
In $M$, for each pair $(p,K)\in\poset\x\N$ and each $K_0\in\N$,
there are $K_1\in\N$ and a sequence $\seq{m_n}{K_0\le n<K_1}$, such that
$X(p,K)\sub\bigcup_{K_0\le n<K_1}U^n_{m_n}$.
\eclm
\bpf
If $X(p,K)\sub\bigcup_{K_0'\le n<K_1}U^n_{m_n}$ for some $K_0'\ge K_0$, then
$X(p,K)\sub\bigcup_{K_0\le n<K_1}U^n_{m_n}$. Thus, we may assume that $K_0\ge K$.
Take $q\le p$ and $K_1$ such that
$q\forces \name{f}(K_0)=K_1$. Extend $q$
so that $K_1$ is in the domain of $q$.
Then
\begin{eqnarray*}
X(p,K) & \sub &
\medset{x\in X}
{p\forces x\in \bigcup_{K_0\le n<\name{f}(K_0)}U^n_{\name{g}(n)}} \sub\\
 & \sub &
\medset{x\in X}
{q\forces x\in \bigcup_{K_0\le n<\name{f}(K_0)}U^n_{\name{g}(n)}} =\\
& = & \bigcup_{K_0\le n<K_1}U^n_{q(n)}.\qedhere
\end{eqnarray*}
\epf

Enumerate $\poset\x\N=\seq{(p_i,N_i)}{i\in\N}$.
Using the claim, pick numbers $K_1$ and $m_n$ for $n<K_1$ such that
$X(p_0,N_0) \sub\bigcup_{n<K_1}U^n_{m_n}$.
Pick numbers $K_2$ and $m_n$ for $K_1\le n<K_2$ such that
$X(p_1,N_1) \sub\bigcup_{K_1\le n<K_2}U^n_{m_n}$.
Pick numbers $K_3$ and $m_n$ for $K_2\le n<K_3$ such that
$X(p_2,N_2) \sub\bigcup_{K_2\le n<K_3}U^n_{m_n}$.
Continuing in this manner, we obtain a sequence $\seq{m_n}{n\in\N}\in M$
in $\N$ such that
$$X=\Un_{i\in\N}X(p_i,K_i)\sub \Un_{n\in\N} U^n_{m_n}.$$

$(3)\Impl (2)$: We will use the following lemma.

\blem\label{lem:game}
Assume that, in the ground model, a space $X$ with points G$_\delta$ satisfies $\sone(\Op,\Op)$.
Assume that $\poset$ is a poset and $\name{\cU}$ is a $\poset$-name for an open cover of $X$, consisting of open sets from the ground model.
For each $p\in\poset$, there are a decreasing sequence
$\seq{q_m}{m\in\N}$ in $\poset$
and a sequence $\seq{U_m}{m\in\N}$ of sets open in the ground model such that:
\be
\item $q_0=p$;
\item $q_{m+1}\forces U_m\in\name{\cU}$ for all $m$; and
\item $\sset{U_m}{m\in\N}$ is a cover of $X$.
\ee
\elem
\bpf
For each condition $q\in\poset$, let
$$\cU_q=\set{U}{\exists r\le q,\ r\forces U\in\name{\cU}}.$$
Then $\cU_q\in M$, and is a cover of $X$.

Define a strategy for ONE in the Rothberger game $\gone(\Op,\Op)$
on $X$, as follows.
Let $q_0=p$.
ONE's first move is the cover $\cU_{q_0}$.
Suppose that TWO responds with an element $U_0\in\cU_{q_0}$.
Then ONE picks, using a fixed choice function on
the nonempty subsets of $\poset$,
a condition $q_1\le q_0$ forcing that $U_0\in \name{\cU}$, and
plays $\cU_{q_1}$. If TWO responds with an element $U_1\in\cU_{q_1}$, then
ONE picks $q_2\le q_1$ forcing that $U_1\in\name{\cU}$,
and plays $\cU_{q_2}$, and so on.

By Pawlikowski's Theorem~\cite[Theorem 1]{Paw},
since $X$ satisfies $\sone(\Op,\Op)$, the strategy thus defined is not a winning strategy.
Let $\seq{q_m}{m\in\N}$ and $\seq{U_m}{m\in\N}$ be the sequences occurring during a play
lost by ONE. Then (1)--(3) hold.
\epf

Let $\poset=\two^{<\alephes}$.
Let $\seq{\name{\cU}_n}{n\in\N}$ be a sequence of $\poset$-names for open
covers of $X$ consisting of ground model open sets.

Fix $n$ and a condition $p\in\poset$. By Lemma~\ref{lem:game},
there are a decreasing sequence $\seq{q_m^{n,p}\in\poset}{m\in\N}$
and a sequence $\seq{U_m^{n,p}}{m\in\N}\in M$ of open subsets of $X$ such that
\be
\item $q_0^{n,p}=p$;
\item $q_{m+1}^{n,p}\forces U_m^{n,p}\in\name{\cU}_n$ for all $m$; and
\item $\set{U_m^{n,p}}{m\in\N}$ is a cover of $X$.
\ee
As $X$ satisfies $\ufin(\Op,\Ga)$, there are for each pair $(n,p)\in \N\x\poset$
a number $k(n,p)$ such that
$$\medset{\Un_{m<k(n,p)}U^{n,p}_m}{(n,p)\in \N\x\poset}\in\Ga(X).$$
By enlarging the numbers $k(n,p)$, we may assume that the displayed enumeration
is bijective.

Let $G$ be $\poset$-generic over $M$.
Fix $n$. The set $\sset{q^{n,p}_{k(n,p)}}{p\in\poset}$ is dense in $\poset$.
Let $p_n$ be a condition such that $q^{n,p_n}_{k(n,p_n)}\in G$.
Then, in $M[G]$, we have that
$$\set{U^{n,p_n}_m}{m<k(n,p_n)}\sub\cU_n.$$
As our enumeration is bijective, we have that
$$\medset{\Un_{m<k(n,p_n)}U^{n,p_n}_m}{n\in\N}\in\Ga(X).$$
This completes the proof.
\epf

\brem
In Theorem~\ref{thm:starpres}, the only implication that uses the premise that the points of
the space are G$_\delta$ is ``$(3)\Impl (2)$''. Since this hypothesis is very mild, we have not
tried to eliminate it.
\erem

Theorem~\ref{thm:starpres} has the following corollary.

\bcor
For spaces with points G$_\delta$, the Gerlits--Nagy property (*)
(equivalently, the conjunction of $\ufin(\Op,\Ga)$ and $\sone(\Op,\Op)$)
is preserved by Cohen forcing.\qed
\ecor


\section{$\Cp$ theory and more applications}\label{sec:cp}

For a space $X$, let $\rmD(X)$ be the
family of all dense subsets of $X$.
Spaces satisfying $\sfin(\rmD,\rmD)$ are also called \emph{selectively separable} or \emph{M-separable},
and spaces satisfying $\sone(\rmD,\rmD)$ are also called \emph{R-separable}---see~\cite{DFamilies} for a summary and references.\footnote{In the paper~\cite{DFamilies},
the family $\rmD$ is defined differently, in order to study additional properties in a uniform manner.
The change in the definition of $\rmD$ does not change the property $\sfin(\rmD,\rmD)$.}
For a space $X$ and a point $x\in X$, let $\Om_x(X)$ be the
family of all sets $A\sub X$ with $x\in \cl{A}\sm A$. A space $X$ has \emph{countable fan tightness}
if $\sfin(\Om_x,\Om_x)$ holds at all points $x\in X$. It has \emph{strong countable fan tightness}
if $\sone(\Om_x,\Om_x)$ holds at all points $x\in X$. When the space $X$ is a topological group,
it suffices to consider $\sfin(\Om_x,\Om_x)$ and $\sone(\Om_x,\Om_x)$ at the neutral element of that
group.

Generalizing results of Scheepers~\cite[Theorems 13 and 35]{coc6},
Bella, Bonanzinga, Matveev and Tkachuk prove in~\cite[Corollary 2.10]{BBMT} that the following
assertions are equivalent for every space $X$ and each $\Sel\in\{\sone,\sfin\}$:
\be
\item $\Cp(X)$ satisfies $\Sel(\rmD,\rmD)$;
\item $\Cp(X)$ is separable and satisfies $\Sel(\Om_0,\Om_0)$;
\item $X$ has a coarser, second countable topology, and satisfies $\Sel(\Om,\Om)$.
\ee

\bcor[CH]\label{cor:1}
There are sets $X,Y\sub\R$ such that the spaces $\Cp(X)$ and $\Cp(Y)$ are \FU{},
and their product $\Cp(X)\x \Cp(Y)$ does not satisfy $\sfin(\rmD,\rmD)$ (or, equivalently,
$\sfin(\Om_0,\Om_0)$).
\ecor
\bpf
By Theorem~\ref{thm:manysol}, there are \gs{}s $X,Y\sub\R$ (so that $\Cp(X)$ and $\Cp(Y)$ are
\FU{}) such that $X\sqcup Y$ does not satisfy $\sfin(\Om,\Om)$, and hence
$\Cp(X)\x \Cp(Y)=\Cp(X\sqcup Y)$ does not satisfy $\sfin(\rmD,\rmD)$~\cite[Theorem 35]{coc6}.
\epf

Corollary~\ref{cor:1} strengthens Babinkostova's Corollary 2.5 in~\cite{Bab09},
where the spaces $\Cp(X)$ and $\Cp(Y)$ satisfy the weaker property $\sone(\rmD,\rmD)$. In fact,
Babinkostova's space are provably not \FU{}. When this extra feature is taken into account, 
the results become incomparable.
Corollary~\ref{cor:1} can be used to reproduce a result of Barman and Dow~\cite[Theorem 2.24]{BarmanDow}.
The Barman--Dow Theorem is identical to Corollary~\ref{cor:2} below, except that their countable
spaces are not topological groups.

\bcor[CH]\label{cor:2}
There are countable abelian \FU{} topological groups $A$ and $B$ such that
the product group $A\x B$ does not satisfy $\sfin(\rmD,\rmD)$ or $\sfin(\Om_0,\Om_0)$.
\ecor
\bpf
It suffices to consider $\sfin(\rmD,\rmD)$. Indeed, according to~\cite[Proposition 2.3(2)]{BBMT},
every separable space with countable fan tightness satisfies $\sfin(\rmD,\rmD)$.

Let $X$ and $Y$ be as in Corollary~\ref{cor:1}. Let $\seq{D_n}{n\in\N}$ be a sequence of countable dense
subsets of $\Cp(X)\x\Cp(Y)$ witnessing the failure of $\sfin(\rmD,\rmD)$ for $\Cp(X)\x\Cp(Y)$.
Let $A$ and $B$ be the groups generated by the projections of $\Un_n D_n$ on the first
and second coordinates, respectively.
As being \FU{} is hereditary, the countable groups $A$ and $B$ are \FU{}.
As $A\x B$ contains $D_0$, it is dense in $\Cp(X)\x\Cp(Y)$.
The sets $D_n$ are contained in $A\x B$, and are dense (in particular) there.
Assume that there are finite sets $F_n\sub D_n$ for $n\in\N$ such that $\Un_nF_n$ is dense in
$A\x B$. Then $\Un_nF_n$ is dense in $\Cp(X)\x\Cp(Y)$, a contradiction.
\epf

To what extent is \CH{} necessary for Theorem~\ref{thm:nonprod}?
Typically, in the field of selection principles, Martin's Axiom suffices to establish
consequences of \CH{}. Surprisingly, this is not the case here.
The following theorem is an immediate consequence of a result of
Barman and Dow~\cite[Theorem 3.3]{BarmanDow2}.
PFA stands for the Proper Forcing Axiom, an axiom that is strictly stronger than Martin's Axiom.

\bthm[PFA]\label{thm:BD}
All finite products of separable metric \gs{}s satisfy $\sfin(\Om,\Om)$.
\ethm
\bpf
According to a result of Barman and Dow~\cite[Theorem 3.3]{BarmanDow2},
PFA implies that all finite products of countable \FU{} spaces satisfy $\sfin(\rmD,\rmD)$.
We consider products of two sets. The generalization to arbitrary finite products
is straightforward.

Assume that $X$ and $Y$ are separable metric \gs{}s and $X\x Y$ does not satisfy $\sfin(\Om,\Om)$.
As the property $\sfin(\Om,\Om)$ is preserved by finite powers~\cite[Theorem 2.5]{coc2},
Lemma~\ref{lem:un2prod} implies that $X\sqcup Y$ does not satisfy $\sfin(\Om,\Om)$. Thus,
by Scheepers's Theorem, the space $\Cp(X)\x\Cp(Y)=\Cp(X\sqcup Y)$ does not satisfy $\sfin(\rmD,\rmD)$.
Continuing as in the proof of Corollary~\ref{cor:2}, we obtain two countable \FU{} spaces whose
product is not $\sfin(\rmD,\rmD)$; a contradiction.
\epf

By the above-mentioned theorem of Bella, Bonanzinga, Matveev and Tkachuk~\cite[Corollary 2.10]{BBMT},
it suffices to assume in Theorem~\ref{thm:BD}, that the \gs{}s have a coarser, second countable
topology.

In the Cohen model, a result stronger than Theorem~\ref{thm:BD}
follows from another result of Barman and Dow~\cite{BarmanDow2}.

\bthm
In the Cohen model, obtained by adding at least $\aleph_2$ Cohen reals to a model of \CH{},
all separable metrizable \gs{}s $X$ have cardinality at most $\aleph_1$.
\ethm
\bpf
Let $X$ be a Tychonoff \gs{}. Then $\Cp(X)$ is \FU{}. Fix a countable dense subset $D$ of $\Cp(X)$.
Then $D$ is \FU{}. According to~\cite[Theorem 3.1]{BarmanDow2}, 
in the Cohen model, all countable \FU{} spaces have $\pi$-weight at most $\aleph_1$.
It follows that the $\pi$-weight of $D$ is at most $\aleph_1$.
By the density of $D$, the $\pi$-weight of $\Cp(X)$ is at most $\aleph_1$.
In a topological group, if $\cU$ is a pseudo-base, then the set $\sset{U\inv\cdot U}{U\in\cU}$
is a local base at the neutral element. Thus, the cardinality of $X$, which is equal to the character of $\Cp(X)$,
is at most $\aleph_1$.
\epf

As $\aleph_1<\fd$ in the Cohen model, the consequence that products of \gs{}s in $\R$ satisfy
$\sfin(\Om,\Om)$ there is trivial, i.e., follows from sheer cardinality considerations.

\medskip

The following theorem solves, in the negative, Problem 3.1 (and thus also Problems~3.2 and 3.3)
of Samet--Tsaban~\cite[{\S}3]{SPMB18}. 
This problem asks whether every set $X\sub\R$ with the Hurewicz property, and with Menger's property in all finite powers, necessarily has the Hurewicz property in all finite powers.
Theorem~\ref{thm:gnoth} also provides a consistently positive answer to Problem~3.4 there,
since adding $\aleph_1$ Cohen reals to a model of \CH{} preserves \CH{}.
A proposed solution of these problems 
in~\cite{SchTall10} is withdrawn in~\cite{STErrata}, for the reasons in
the discussion following Problem~\ref{prb:weiss}.

\bthm\label{thm:gnoth}
In any model obtained by adding uncountably many Cohen reals to a model of \CH{},
there is a set $X\sub\R$ such that $X$ satisfies $\sone(\Tau,\Ga)$ and $\sone(\Om,\Om)$,
but its square $X^2$ does not satisfy $\ufin(\Op,\Ga)$.
\ethm
\bpf
In the ground model, using \CH{}, let $X$ be the set in the proof of Theorem~\ref{thm:manysol}.
Move to the generic extension.
By Theorem~\ref{thm:gpres}, the set $X$ remains the union of two \gs{}s.
Thus, $X$ satisfies $\sone(\Tau,\Ga)$.
All finite powers of ground model sets, including $X$, satisfy $\sone(\Op,\Op)$ in the
extension~\cite[Theorem 11]{SchTall10}.
Equivalently, $X$ satisfies $\sone(\Om,\Om)$.
By Theorem~\ref{thm:manysol}, in the ground model, the square $X^2$ does not satisfy $\sfin(\Op,\Op)$,
and thus does not satisfy $\ufin(\Op,\Ga)$. It follows that, in the extension, the square $X^2$ does not satisfy 
$\ufin(\Op,\Ga)$~\cite[Theorem 37]{SchTall10}.
\epf

Similarly, we have the following.

\bthm
In any model obtained by adding
uncountably many Cohen reals to a model of \CH{}, there are \gs{}s $X,Y\sub\R$
such that $X\x Y$ satisfies $\sone(\Om,\Om)$ but not $\ufin(\Op,\Ga)$.\qed
\ethm

For a space $X$, let $\cD\in\rmD_\Ga(X)$ if $\cD$ is infinite, and for each open set $U$ in $X$,
$U$ intersects all but finitely many members of $\cD$. Spaces satisfying $\sfin(\rmD,\rmD_\Ga)$
are also called \emph{H-separable} (e.g.,~\cite{BBM}). Also, for $x\in X$,
let $\Ga_x$ be the family of all countable sets converging to $x$.
Spaces satisfying $\sone(\Ga_x,\Ga_x)$ are also called \emph{$\alpha_2$ spaces}.

\bcor\label{cor:co}
In any model obtained by adding uncountably many Cohen reals to a model of \CH{},
there is a set $X\sub\R$ such that the space $\Cp(X)$
satisfies $\sone(\rmD,\rmD)$ and $\sone(\Ga_0,\Ga_0)$, but not $\sfin(\rmD,\rmD_\Ga)$.
\ecor
\bpf
Let $X$ be the set from Theorem~\ref{thm:gnoth}.
As $X$ satisfies $\sone(\Om,\Om)$, the space $\Cp(X)$ satisfies $\sone(\rmD,\rmD)$
\cite[Theorem 13]{coc6}.
As $X$ satisfies $\sone(\Ga,\Ga)$, the space $\Cp(X)$ satisfies $\sone(\Ga_0,\Ga_0)$
\cite[Theorem 4]{CpHure}.
As $X^2$ does not satisfy $\ufin(\Op,\Ga)$, the space $\Cp(X)$ does not satisfy $\sfin(\rmD,\rmD_\Ga)$
\cite[Theorem 40]{BBM}.
\epf

By the usual method used in the earlier proofs, Corollary~\ref{cor:co} has the following consequence.

\bcor
In any model obtained by adding uncountably many Cohen reals to a model of \CH{},
there is a countable abelian topological group $A$ satisfying $\sone(\rmD,\rmD)$ and $\sone(\Ga_0,\Ga_0)$,
but not $\sfin(\rmD,\rmD_\Ga)$. \qed
\ecor


\section{The product of an unbounded tower set and a Sierpi\'nski set}

We conclude this paper with a proof that, for each unbounded tower
$T=\set{x_\alpha}{\alpha<\fb}\sub\roth$
and each Sierpi\'nski set $S$, the product space $(T\cup\Fin)\x S$ satisfies $\sone(\Ga,\Ga)$.
In fact, we prove a more general result.

For each unbounded tower $T=\set{x_\alpha}{\alpha<\fb}\sub\roth$, the set
$T\cup\Fin$ satisfies $\sone(\Ga,\Ga)$
(implicitly in~\cite[Theorem 6]{alpha_i}, and explicitly in~\cite[Proposition 2.5]{CBC}).
The existence of unbounded towers of cardinality $\fb$ follows from the
existence of unbounded towers of any cardinality~\cite[Proposition 2.4]{BBC}.
Examples of hypotheses implying
the existence of unbounded towers are $\ft=\fb$ or $\fb<\fd$~\cite[Lemma 2.2]{BBC}.

The property $\sone(\BG,\BG)$ is equivalent to the Hurewicz property for countable Borel
covers, and also to the property that all Borel images in the Baire space $\NN$ are bounded
\cite[Theorem 1]{CBC}.

\bthm\label{thm:tgg}
Let $T=\set{x_\alpha}{\alpha<\fb}\sub\roth$ be an unbounded tower.
For every space $Y$ satisfying $\sone(\BG,\BG)$,
the product space
$(T\cup\Fin)\x Y$ satisfies $\sone(\Gamma,\Gamma)$.
\ethm
\bpf
Let $\cU=\set{U_n}{n\in\N}\in\Ga((T\cup\Fin)\x Y)$.
For a finite set $s\sub\N$ and $n\in\N$, let
$$[s,n]=\set{x\sub\N}{x\cap \{0,\ldots,n-1\}=s}\cap (T\cup\Fin).$$
By shrinking the elements of $\cU$, we may assume that
$U_n\cap (\{n\}\x Y)=\emptyset$ for all $n$.
Consider the functions $f{},g{}\colon Y\to\NN$, defined by
\begin{eqnarray*}
f{}(y)(n) & = & \max\set{k}{P(\{0,\ldots,k-1\})\x\{y\}\sub U_n},\\
g{}(y)(n) & = & \min\set{l\ge n}{\forall s\in P(\{0,\ldots,f{}(y)-1\}),\ [s,l]\x\{y\}\sub U_n}.
\end{eqnarray*}
By our assumption on $\cU$, we have that $f{}(y)(n)\le n$.
As $\cU\in\Ga((T\cup\Fin)\x Y)$, the sequence $\seq{f{}(y)(n)}{n\in\N}$
converges to infinity for each $y\in Y$.

\bclm\label{borels}
The function $f{}$ is Borel, and there is a Borel function
$h{}\colon Y\to\NN$ such that $g{}(y)(n)\le h{}(y)(n)$ for all $y\in
Y$ and all $n$.
\eclm
\bpf
The function $f{}$ is Borel,
since the preimages under $f{}$ of the standard basic open subsets of $\NN$ are finite
intersections of subsets of $Y$ which are either closed or open.

Represent each open set $U_n$ as an increasing union $\bigcup_{k}U_{n,k}$ of clopen sets.
Let $\overline{\N}$ be the set $\N\cup\{\infty\}$, with the discrete topology.
Define a function
$\Phi\colon Y\to \bigl(\overline{\N}^\N\bigr)^\N$ as follows: $\Phi(y)(n)(k)=\infty$
if
$P(\{0,\ldots,f{}(y)(n)-1\})\x\{y\}\not\sub U_{n,k}$, and if not, then
$\Phi(y)(n)(k)$ is the minimal $l$ such that
$[s,l]\x\{y\}\sub U_{n,k}$ for all $s\sub \{0,\ldots,f{}(y)(n)-1\}$.
Since $P(\{0,\ldots,f{}(y)(n)-1\})\x\{y\}\sub U_n$, by the definition
of $f{}$, there is $k$ such that $P(\{0,\ldots,f{}(y)(n)-1\})\x\{y\}\sub U_{n,k}$.
Thus, the set $\sset{k}{\Phi(y)(n)(k)=\infty}$ is finite.
Moreover, the sequence $\seq{\Phi(y)(n)(k)}{k\in\N}$ is nonincreasing (we assume that
$i<\infty$ for all $i$), and $\Phi(y)(n)(k)\ge g{}(y)(n)$ for all $k$.
Set $h{}(y)(n)=\min\set{\Phi(y)(n)(k)}{k\in\N}$. It follows that
$h{}(y)(n)\ge g{}(y)(n)$ for all $n$. Thus, it suffices to
prove that $h{}\colon Y\to\NN$ is Borel, which follows as soon as we prove that
$\Phi\colon Y\to \bigl(\overline{\N}^\N\bigr)^\N$ is Borel.

Fix $n,k\in\N$ and $m\in\overline{\N}$. We need to show
hat the set $A=\set{y\in Y}{\Phi(y)(n)(k)=m}$ is Borel.
Consider the two possible cases.

\noindent\emph{Case 1: $m=\infty.$}
In this case,
\begin{eqnarray*}
A & = & \set{y}{P(\{0,\ldots,f{}(y)(n)-1\})\x\{y\}\not\sub U_{n,k}}=\\
& = &
\bigcup_{l<n}\big(\set{y\in Y}{f{}(y)(n)=l}\cap\set{y\in Y}{P(\{0,\ldots,l-1\})\x\{y\}\not\sub U_{n,k}}\big)=\\
& = & \bigcup_{l<n}\big(\set{y\in Y}{f{}(y)(n)=l}\cap\bigcup_{s\sub \{0,\ldots,l-1\}}\set{y\in Y}{(s,y)\notin U_{n,k}}\big).
\end{eqnarray*}
As the function $f{}$ is Borel, the set $\sset{y\in Y}{f{}(y)(n)=l}$ is Borel.
The set $\sset{y\in Y}{(s,y)\notin U_{n,k}}$ is a clopen subset of $Y$ for all $s\sub l$.
Thus, $A$ is Borel.

\noindent\emph{Case 2: $m\in\N.$} In this case,
\begin{eqnarray*}
A & = & \set{y\in Y}{\forall s\sub \{0,\ldots,f{}(y)(n)-1\},\ ([s,m]\x\{y\}\sub U_{n,k})}\cap \\
& & \cap
\set{y\in Y}{\exists s\sub \{0,\ldots,f{}(y)(n)-1\},\ ([s,m-1]\x\{y\} \not\sub U_{n,k}) }=\\
& = &
\bigcup_{l<n} \big(\set{y\in Y}{\forall s\sub \{0,\ldots,l-1\},\ ([s,m]\x\{y\}\sub U_{n,k})}\cap\\
& & \cap \set{y\in Y}{f{}(y)(n)=l}\big)\cap\\
& & \cap \bigcup_{l<n}\big( \set{y\in Y}{\exists s\sub \{0,\ldots,l-1\},\ ([s,m-1]\x\{y\} \not\sub U_{n,k}) }\cap \\
& & \cap\set{y\in Y}{f{}(y)(n)=l} \big).
\end{eqnarray*}
As the function $f{}$ is Borel, the latter set is Borel.
Indeed, for each $V\sub T\cup \Fin$, the set
$\set{y\in Y}{V\x\{y\}\sub U}$ is closed whenever $U\sub (T\cup \Fin)\x Y$ is closed.
\epf

\bclm\label{bounds}
There is an increasing function $c\in\NN$
such that, for each $y\in Y$,
$$c(n)\le f{}(y)(c(n+1))\le h{}(y)(c(n+1))<c(n+2)$$
for all but finitely many $n$.
\eclm
\bpf
Consider the map $f'{}\colon Y\to\NN$, defined by
$f'{}(y)(n)=\min\set{f{}(y)(l)}{l\ge n}$.
Then the set $f'(Y)\sub\NN$ consists of nondecreasing
unbounded sequences.
Set
$$f''{}(y)(k)=\min\set{n}{f'{}(y)(n)\ge k}.$$
Then $f''{}\colon Y\to\NN$ is a Borel map, and hence $f''{}(Y) $
is bounded by some increasing function $a'\in\NN$.
Let $a(n)=\min\set{k}{a'(k)\ge n}$.
Then $a\le^* f'{}(y)\le^* f{}(y)$ for all $y\in Y$.

Since $Y$ satisfies $\sone(\BG, \BG)$ and $h{}$ is
Borel, there is an increasing $b\in\NN$ such that
$h{}(y)\le^* b$ for all $y\in Y$.
Let $c(0)=1$, and
$$c(n+1)=\max\{\min\set{l}{a(l)\ge c(n)},b(c(n))\}+1.$$
We claim that $c$ is as required. Indeed, fix $y\in Y$ and find $n$
such that $a(m)\le h{}(y)(m)\le g{}(y)(m)\le b(m)$ for all
$m\ge n$. For $m\ge n$, as
$c(m+1)\ge\min\set{l}{a(l)\ge c(m)}$
and
$a$ is nondecreasing,
we have that
$f{}(y)(c(m+1))\ge a(c(m+1))\ge c(m)$,
and the inequality $h{}(y)(c(m+1))\le
b(c(m+1))<c(m+2)$ follows.
\epf

Let $\seq{\cU_k}{k\in\N}$ be a sequence in $\Ga((T\cup\Fin)\x Y)$,
where $\cU_k=\seq{U^k_n}{n\in\N}$ for all $k$.

\bclm\label{lem:small_mistake_allowed}
Suppose that for every sequence $\seq{\cV_k}{k\in\N}$
in $\Ga((T\cup\Fin)\x Y)$,
where $\cV_k=\seq{V^k_n}{n\in\N}$ for all $k$, there exists a sequence
$\seq{n_k}{k\in\N}$ in $\N$ such that $\seq{V^k_{n_k}}{k\in\N}\in
\Gamma(A\x Y)$ for some $A$ containing $\Fin$ with $|T\sm A|<\fb$. Then
$(T\cup\Fin)\x Y$ is $S_1(\Gamma,\Gamma)$.
\eclm
\bpf
First let us note that the following statement may be obtained simply by splitting
each $\cV_n$ into countably many disjoint infinite pieces
and applying the assumption to the sequence of pieces:\\
for every sequence $\seq{\cV_k}{k\in\N}$ in $\Ga((T\cup\Fin)\x Y)$
there exists a sequence
$\seq{\cV'_k}{k\in\N}$ such that $\cV'_k$ is an infinite subset
of $\cV_k$ and
 $\bigcup_{k}\cV'_k\in
\Gamma(A\x Y)$ for some $A$ containing $\Fin$ with $|T\sm A|<\fb$.

Fix $\alpha_0<\mathfrak b$ and a sequence $\seq{\cV_k}{k\in\N}$ in $\Ga((T\cup\Fin)\x Y)$. Since the set
$\set{x_\xi}{\xi<\alpha_0}\times Y$ is $S_1(\mathrm{B}_\Gamma,\mathrm{B}_\Gamma)$,
there exists a sequence
$\seq{\cW^0_k}{k\in\N}$ such that $\cW^0_k$ is an infinite subset of
$\cV_k$ and
 $\bigcup\set{\cW^0_k}{k\in\mathbb N}\in\Gamma(\set{x_\xi}{\xi<\alpha_0}\times Y)$.
Applying (the reformulation of) our assumption to the sequence $\seq{\cW^0_k}{k\in\N}$
in $\Ga\bigl((T\cup\Fin)\x Y\bigr)$, we can find a sequence
$\seq{\cV^0_k}{k\in\N}$ such that
$\cV^0_k$ is an infinite subset
of $\cW^0_k$ and
 $\bigcup_{k}\cV^0_k\in
\Gamma(A\x Y)$ for some $A$ containing $\Fin$ with $|T\sm A|<\fb$.
It follows that $\bigcup_{k}\cV^0_k\in
\Gamma((\set{x_\xi}{\xi<\alpha_0}\cup \set{x_\xi}{\xi>\alpha_1}\cup\Fin)\x Y)$
for some $\alpha_1>\alpha_0$.

Applying the same argument infinitely many times we can get
an increasing sequence $\seq{\alpha_n}{n\in\mathbb N}$ of ordinals
below $\mathfrak b$,
and for every $n$ a sequence
$\seq{\cV^n_k}{k\in\N}$ such that
$\cV^n_k$ is an infinite subset
of $\cV^{n-1}_k$ and
 $\bigcup_{k}\cV^n_k\in
\Gamma((\set{x_\xi}{\xi<\alpha_n}\cup \set{x_\xi}{\xi>\alpha_{n+1}}\cup\Fin)\x Y)$.
Let us select $V_k\in \cV^k_k\sm\{V_0,\dots,V_{k-1}\}$ for all $k$.
Then $V_k\in \cV_k$ and
$\set{V_k}{k\in\N}$ is easily seen to be in
$\Gamma((T\cup\Fin)\x Y)$.
\epf

By Lemma~\ref{lem:small_mistake_allowed},
it suffices to find a sequence $\seq{n_k}{k\in\N}$ in $\N$
such that $\seq{U^k_{n_k}}{k\in\N}\in
\Gamma(A\x Y)$ for some $A$ containing $\Fin$ with $|T\sm A|<\fb$.

For each $k\in\N$, let $c_k\in\NN$ be such as in Claim~~\ref{bounds},
where $\cU$ is replaced with $\cU_k$, and let $f_k$ and $h_k$ be the associated functions.
Consider the function
$d\colon Y\to\NN$ defined by
$$ d(y)(k)=\min\set{n}{\forall m\ge n,\ \big(c_k(m)\le f_{{}k}(y)(c_k(m+1))<h_{{}k}(y)(c_k(m+1))<c_k(m+2)\big) }. $$
Since the functions $f_{{}k}$ and $h_{{}k}$ are Borel, so is the function $d$, and hence
there is an increasing $x\in\NN$ such that $d(y)\le^* x$ for
all $y\in Y$. We may assume that
$c_{k+1}(x(k+1))>c_k(x(k)+2)$ for all $k$.
Let $\alpha<\fb$ be such that the set
$I=\set{k}{x_\alpha\cap [c_k(x(k)), c_k(x(k)+2))=\emptyset}$
is infinite. Fix $\beta\ge\alpha$ and $y\in Y$,
and find $k_0$ such that
$x_\beta\sm x_\alpha\sub k_0$ and $d(y)(k)\le x(k)$ for all $k\ge k_0$.
Then, for all $k\ge k_0$ in $I$, we have that
$x_\beta\cap [c_k(x(k)), c_k(x(k)+2))\sub x_\alpha\cap [c_k(x(k)), c_k(x(k)+2))=\emptyset$.
Consequently,
$x_\beta\cap [f_{{}k}(y)(c_k(x(k)+1)),h_{{}k}(y)(c_k(x(k)+1)))=\emptyset$,
and hence $x_\beta\cap [f_{{}k}(y)(c_k(x(k)+1)), g_{{}k}(y)(c_k(x(k)+1)))=\emptyset$.
Thus,
$$(x_\beta,y)\in [x_\beta\cap f_{{}k}(y)(c_k(x(k)+1)), g_{{}k}(y)(c_k(x(k)+1))].$$
By the definitions of $f_{{}k}$ and $g_{{}k}$, the latter open set is a subset of
$U^k_{c_k(x(k)+1)}$. Therefore, for every
$\beta\ge\alpha$ and $y\in Y$, we have that $(x_\beta,y)\in U^k_{c_k(x(k)+1)}$
for all but finitely many $k\in I$.
As the covers $\cU_k$ get finer with $k$, this completes our proof.
\epf

As the unbounded set $T$ in Theorem~\ref{thm:tgg} is a Borel subset of
the space $T\cup\Fin$, the latter space does not satisfy $\sone(\BG,\BG)$. In particular,
it is not productively $\sone(\BG,\BG)$.

\bprb
Let $T=\set{x_\alpha}{\alpha<\fb}$ be an unbounded tower.
Is the space $T\cup\Fin$, provably, productively $\sone(\Ga,\Ga)$?
Is this the case assuming \CH{}?
\eprb

\ed

\end{document}